\documentstyle{article}
\input{amssym.def}
\input{amssym.tex}
\newtheorem{lemma}{Lemma}[section]
\newtheorem{proposition}[lemma]{Proposition}
\newtheorem{theorem}[lemma]{Theorem}
\newtheorem{corollary}[lemma]{Corollary}
\newtheorem{definition}[lemma]{Definition}
\newcommand{\freeline}{\vspace{0.5cm}}
\newcommand{\freespace}{\hspace{0.3cm}}
\newcommand{\surjection}{\mbox{ $\rightarrow\!\!\!\!\!\!\!\rightarrow$ }}
\newcommand{\QED}{\mbox{\hspace*{\fill}$\Box$}}
\newcommand{\power}[2]{\mbox{$#1\/[[\/#2\/]]$}}
\newcommand{\poly}[2]{\mbox{$#1\/[\/#2\/]$}}

\newcommand{\multpower}[3]{\mbox{$#1\/[[\/#2_{1}, \ldots ,#2_{#3} \/]]$}}
\newcommand{\Opower}[2]{\mbox{${\cal O}_{#1}\/[[\/#2\/]]$}}
\newcommand{\spower}[3]{\mbox{$\power{#1}{#2}^{\oplus #3}$}}
\newcommand{\negpower}[2]{\mbox{$#1\/((\/#2\/))$}}
\newcommand{\snegpower}[3]{\mbox{$\negpower{#1}{#2}^{\oplus #3}$}}
\newcommand{\negpoly}[2]{\mbox{$#1\/(\/#2\/)$}}
\newcommand{\Onegpower}[2]{\mbox{${\cal O}_{#1}\/((\/#2\/))$}}

\title{On a relative version of the Krichever correspondence}
\author{Ines Quandt}
\date{February 18, 1997}
\begin{document}
\maketitle

\section*{Preface}
This PhD thesis is the result of my work in  the Graduiertenkolleg 
"Geometrie und Nichtlineare Analysis" at  
Humboldt University Berlin and in the DFG project KU 770/1-3. 
\freeline\\
It is published in the {\em Bayreuther Mathematische Schriften} {\bf 52} 
(1997), p.1-74.
\freeline\\
At this point, I would like to express my thanks to all of the people who 
supported my mathematical development. 
\vspace{0.5cm}\\
My special thanks go to Doz.~Dr.~sc.~W.~Kleinert, my thesis advisor and my 
professor since the very 
beginning of my studies.
He also introduced me to the fascinating area of algebraic geometry, turned 
my attention to the theory of 
evolution equations, and  kindly supported my work on this thesis.
\vspace{0.5cm}\\
The subject of the present work has been suggested to me by Prof.~M.~Mulase, whom 
I would like to express my  gratitude for his interest in my work and 
the inspiration and the encouragement he gave me.
\vspace{0.5cm}\\
I gratefully thank Prof.~H.~Kurke for his keen interest in my work and lots 
of valuable hints.
The discussions with him have been a wonderful help during the completion of 
this paper.
\vspace{0.5cm}\\
My special thanks also go to  G.~Hein , A.~Matuschke, Dr.~M.~Pflaum and
D.~Ro\ss berg for numerous inspiring
 discussions.
\vspace{1.5cm}\\
Berlin, October 1996 \hfill Ines Quandt
\newpage
\tableofcontents
\newpage
\addtocounter{section}{-1}
\section{Introduction}
The aim of this paper is to construct a link ranging from a class of  sheaves 
on curves over some base scheme via infinite Grassmannians to 
commutative algebras of differential operators and evolution equations.
\freeline

The idea of studying relationships between algebraic curves, algebras 
of differential operators and  partial differential equations is not 
new. This connection has been studied already at the beginning 
of the 
20th century by G.~Wallenberg, I.~Schur, J.~L.~Burchnall 
and 
T.~W.~Chaundy. G.~Wallenberg \cite{W} tried to find all commuting pairs of ordinary 
differential operators.
During his classification of commuting operators $P$ and $Q$ of order 2 and 3, 
respectively, he found a 
certain relation to a plane cubic curve. However, Wallenberg did not 
continue to explore this relation.

The motivation to study commutative algebras instead of commuting pairs of 
differential operators was given by
 I.~Schur \cite{Sc} in 1905, when he proved the following remarkable fact:
\freeline\\
{\em Let $P$ be an ordinary differential operator of order greater than zero, 
and let $B_{P}$ be the set of all
 differential operators which commute with $P$. Then $B_{P}$ is a commutative algebra. }
\freeline

The work of J.~L.~Burchnall 
and 
T.~W.~Chaundy (\cite{BC1}-\cite{BC3}) about the relations between commuting 
differential operators and affine
 algebraic curves is extensive.
For example, they proved that for commuting $P$ and $Q$ of positive order, the 
ring $\poly{\Bbb{C}}{P,Q}$ has 
dimension 1.
Furthermore, they analyzed at length certain examples of affine curves and the 
related differential operators. 
It is remarkable that a large part of the  methods which have been systematically 
developed  more than 50 years
 later, virtually already exist in these early papers by Burchnall and Chaundy, 
although mainly in examples.
\freeline

After a long period of stagnation, another break-through came with the 
work of  P.~Lax \cite{L} about isospectral deformations of differential operators
 in the late 60's. 

It has been realized, already in the early stage of the theory, that
 a commutative algebra of differential operators carries a lot more
 information than only its algebraic-geometric spectrum. 

In the 70's, I.~M.~Krichever analyzed the behaviour of an operator at 
infinity, i.e., he constructed the line bundle on a complete curve corresponding 
to a given algebra by special extensions of the trivial bundle to the
 point at infinity (see \cite{K}, \cite{K1}, \cite{K2}). His approach may be
 considered as the source of the algebraic-geometric correspondences
 established later on. 

Almost simultaneously, and inspired by the work of Krichever \cite{K},
 D.~Mumford \cite{Mum1} established a correspondence between pointed curves equipped
 with a fixed line bundle and certain commutative algebras of 
differential operators. This article and  the 
one of J.-L.~Verdier \cite{V} also contain the very first constructions 
in the case of higher rank vector bundles.
\freeline

Later on, infinite Grassmannians emerged in the study of evolution 
equations and their relations to  vector bundles on curves. Without 
claiming to be complete, let us just mention the work of M.~Sato 
\cite{S}, E.~Previato, G.~Segal and G.~Wilson  (\cite{PW}, \cite{SW}) and 
M.~Mulase \cite{M2}. Among one of the culmination points of this newly 
established theory was  the complete classification of elliptic 
commutative algebras of ordinary differential operators, and one 
of the affirmative solutions of the Schottky problem by M.~Mulase 
in \cite{M1}.\freeline

In  the approach of M.~Mulase, some  questions arise quite naturally: 
Can one generalize the correspondence between vector bundles on curves
 and elements of infinite Grassmannians to the case where the curve is 
not defined over a field, but for example over a $k$-algebra? In this 
case, does one also get a correspondence with commutative algebras of 
differential operators with coefficients in more general rings? The present paper 
gives an affirmative answer to both questions.

These 
questions are interesting from two points of view. First, the 
established correspondence enables us to construct certain classes 
of commutative algebras of partial differential operators. On the other 
hand, it gives us a powerful tool for the study of degenerations of
 curves and vector bundles in terms of differential operators and 
differential equations.
\freeline

In the recently published paper \cite{AMP}
 the authors also give 
a generalization of infinite Grassmannians to the relative case, which 
overlaps with ours in the case where the base scheme is defined over a 
field, and they investigate this Grassmannian from the point of view 
of representation theory.
\freeline

The method presented in our paper serves the purpose of generalizing 
the techniques developed in \cite{M1} to the relative case. Therefore, 
this article can be considered as a general reference and often will  not be 
quoted  in the sequel.
\freeline

A preliminary version of this paper appeared as a preprint \cite{Q}.
In addition to some slight corrections, a great extension has been made:
In the preprint only integral noetherian base schemes were allowed, whereas 
now we only need to assume the base 
scheme to be locally noetherian. 
\freeline

The paper is organized as follows: 
 The first chapter  is devoted to the generalization of the notion of
 infinite Grassmannians. In the second and  third, the link with sheaves 
on relative curves is established. The fourth  chapter  aims
 at illustrating this correspondence. In chapter 5 we give a complete 
characterization of commutative elliptic algebras of differential operators 
with coefficients in the ring of formal power series over a $k$-algebra,
 $k$ being a  field of characteristic zero. The 
appendix is included in order to help those readers who are interested 
in the details of the computations.
\newpage

\section{Relative infinite Grassmannians}

To begin with, we want to generalize the notion of infinite Grassmannians. This 
can be done for arbitrary base schemes.
To this end,  let $S$ be any algebraic scheme.
We denote by $\Opower{S}{z}$ the  sheaf defined by 
$$\Opower{S}{z}(U) := \power{{\cal O}_{S}(U)}{z},$$ 
and $\Onegpower{S}{z}$ is defined as the sheaf
$$\Onegpower{S}{z}(U) := \negpower{{\cal O}_{S}(U)}{z},$$ 
where $\negpower{{\cal O}_{S}(U)}{z}$ stands for the ring of formal 
Laurent series in $z$ with coefficients in ${\cal O}_{S}(U)$. 
$\Onegpower{S}{z}$ has a natural filtration by subsheaves of the 
form $\Opower{S}{z}\cdot z^{n}$. 
\freeline\\
{\bf Remark}\freespace
Neither $\Opower{S}{z}$ nor $\Onegpower{S}{z}$ are quasicoherent sheaves of 
${\cal O}_{S}$-modules.
To see this, take an open affine subset $Spec(R)$ of $S$ and choose an element 
$f\in R$ which is not a zero divisor.
Then there is a natural inclusion
$$\power{R}{z}_{f} \hookrightarrow \power{R_{f}}{z}$$
which, in general, is not an isomorphism. Therefore, $\Opower{S}{z}$ is not quasicoherent.
The same holds true for $\Onegpower{S}{z}$.

However, for each integer $n\in\Bbb{Z}$, the quotient sheaf
$$\Onegpower{S}{z}/\Opower{S}{z}\cdot z^{n} \cong 
\bigoplus_{m<n} {\cal O}_{S}\cdot z^{m}$$
is quasicoherent.

\begin{definition}
Let $U$ be an open subset of $S$ and $v$ a local section of 
$\Onegpower{S}{z}^{\oplus r}$ over $U$, for some $r$. Then the 
{\em order of $v$} is defined to be the minimum integer $n$
 such that $v \in \power{{\cal O}_{S}(U)}{z}^{\oplus r}
\cdot z^{-n}$. 

If $\cal V$ is a subsheaf of $\Onegpower{S}{z}^{\oplus r}$, 
then we define
$${\cal V}^{(n)} := {\cal V} \cap \Opower{S}{z}^{\oplus r}
\cdot z^{-n}.$$
\end{definition}
%
%
%Remark
%
%
{\bf Remark}
$\Onegpower{S}{z}$ acts on $\Onegpower{S}{z}^{\oplus r}$ 
by the natural assignment \begin{equation}
\begin{array}{cccr}
\negpower{{\cal O}_{S}(U)}{z} \times 
\negpower{{\cal O}_{S}(U)}{z}^{\oplus r} & \longrightarrow &
\negpower{{\cal O}_{S}(U)}{z}^{\oplus r} &{}\\
(f,g=(g_{1},\ldots,g_{r})) & \longmapsto & (f\cdot g_{1},\ldots,
f\cdot g_{r}) &.
\end{array}
\end{equation}
For the  multiplication defined this way we have the estimate: 
$$ ord(f\cdot g) \leq ord(f) + ord(g).$$
If $S$ is an integral scheme, then both sides are equal.

For more properties of power series and Laurent series with coefficients in 
arbitrary rings, the reader is referred to  the appendix.
%
%
%Grassmannian
%
%
\begin{definition}
The {\em Grothendieck group} $K(S)$ is defined to be the quotient of the free 
abelian group generated by all coherent 
sheaves on $S$, by the subgroup generated by all expressions 
$${\cal F} - {\cal F}' -{\cal F}''$$
whenever there is an exact sequence 
$$0 \rightarrow {\cal F}' \rightarrow {\cal F} \rightarrow {\cal F}''\rightarrow 0$$
of coherent sheaves on $S$. If $\cal F$ is a coherent sheaf on $S$ then we denote by 
$\gamma({\cal F})$ its image in $K(S)$.
\end{definition}

\begin{definition}
For any natural number $r$, integer $\alpha$ and element $F\in K(S)$, we 
define {\em the infinite Grassmannian of rank $r$, index $F$ and 
level $\alpha$ over $S$} to be the set ${\frak G}^{r}_{F,\alpha}(S)$
 consisting of all quasicoherent subsheaves of ${\cal O}_{S}$-modules
${\cal W} \subseteq \snegpower{{\cal O}_{S}}{z}{r}$ such that 
${\cal W} \cap \spower{{\cal O}_{S}}{z}{r}\cdot z^{-\alpha}$ 
and the quotient 
$\snegpower{{\cal O}_{S}}{z}{r}/({\cal W}+
\spower{{\cal O}_{S}}{z}{r}\cdot z^{-\alpha}) $
are coherent and furthermore:
$$F = 
\gamma({\cal W} \cap \spower{{\cal O}_{S}}{z}{r}\cdot z^{-\alpha})
 - 
\gamma(\snegpower{{\cal O}_{S}}{z}{r}/({\cal W}+
\spower{{\cal O}_{S}}{z}{r}\cdot z^{-\alpha})).$$
\end{definition}
{\bf Remark}\freespace
The introduction of the level has a merely technical meaning. It will be used only in 
Chapter \ref{Fam. DO}. \freeline

Now we introduce the concept of a Schur pair.
%
%
%Schur pair
%
%
\begin{definition}
By a {\em Schur pair of rank $r$ and index $F$ over  $S$} we mean a pair $({\cal A},{\cal W})$ 
consisting of elements ${\cal A} \in {\frak G}^{1}_{G}(S)$, for 
some $G\in K(S)$, and ${\cal W} \in {\frak G}^{r}_{F}(S)$ such 
that 
\begin{itemize}
\item ${\cal A}$ is a sheaf of ${\cal O}_{S}$ - subalgebras of 
$\Onegpower{S}{z}$,
\item The natural action of \negpower{{\cal O}_{S}}{z} on 
\snegpower{{\cal O}_{S}}{z}{r} induces an action of ${\cal A}$ 
on ${\cal W}$, i.e., 
${\cal A}\cdot{\cal W}\subseteq{\cal W}$.
\end{itemize}
We denote by $\frak{S}^{r}_{F}(S)$ the set of Schur 
pairs of rank $r$ and index $F$ over $S$.
\end{definition}
{\bf Remark 1}\freespace
Let us include here a remark on Grothendieck groups. First assume that $S$ is integral. 
Then there is a 
surjective group homomorphism 
$$rk: K(S) \longrightarrow \Bbb{Z}$$
induced by the map
$$\gamma({\cal F}) \longmapsto rk({\cal F}),$$
where $rk({\cal F})$ denotes the rank of $\cal F$ at the generic point of $S$.  If  
$S$ equals $Spec(k)$,
 for some field $k$, then the  homomorphism $rk$ is an isomorphism.

If $S$ is  reduced we still can define a ``multirank'' by taking the rank at every 
irreducible component.
 However, if $S$ is not reduced, the  rank is no longer well-defined. That is why we use 
the Grothendieck
  group to define Grassmannians.
\newpage
{\bf Remark 2}\label{Grass Mulase}\freespace
How are the notions of infinite Grassmannians and Schur 
pairs related to those introduced by M.~Mulase? 
To answer this question, consider the embedding 
$$
\begin{array}{rcl}
\Onegpower{S}{z} &\hookrightarrow &\Onegpower{S}{y}\\
z & \mapsto & y^{r}.
\end{array}
$$
This leads to the  natural identification:
$$
\begin{array}{rcl}
\Onegpower{S}{z}^{\oplus r} & = & \bigoplus_{i=0}^{r-1} 
\Onegpower{S}{y^{r}}\cdot y^{i}\\
&=& \Onegpower{S}{y}.
\end{array}
$$

In particular, let $S=Spec(k)$ for some field $k$. Then 
${\frak G}^{r}_{F}(S)$ consists of all subspaces $W\subset 
\negpower{k}{z}^{\oplus r}$ such that the composition of morphisms 
$$W\hookrightarrow \negpower{k}{z}^{\oplus r} \surjection 
\negpower{k}{z}^{\oplus r}/ \power{k}{z}^{\oplus r}$$
is a Fredholm operator of index $rk(F)$. In view of the previous 
identifications, this redefines the notion of an infinite 
Grassmannian as used by M.~Mulase (cf. \cite{M1}). 

However,   comparing the notions of Schur pairs,  we see
 that our definition is a restricted version of the one given by Mulase. 
Later we will see that in the case where the ground field is of characteristic 
zero, this restriction is not substantial. For more details see 
Section \ref{Fam. DO}. 
\freeline

\begin{definition}
Let $\alpha: S\rightarrow S'$ be a morphism, and $({\cal A}',
{\cal W}')$
 a Schur pair of rank $r$ over $S'$. We denote by $\alpha^{(*)}{\cal A}'$ 
(resp. $\alpha^{(*)}{\cal W}'$) the image of ${\cal A}'$ (resp. ${\cal W}'$) under the map
$$\alpha^{(*)} : \Onegpower{S'}{z}^{\oplus r} \rightarrow 
\Onegpower{S}{z}^{\oplus r}$$
which is given by the pull-back of the coefficients.
\end{definition}
{\bf Remark}\freespace
$(\alpha^{(*)}{\cal A}',\alpha^{(*)}{\cal W}')$ is a Schur pair
 of rank $r$ over $S$.
\freeline

Now we can define  homomorphisms of Schur pairs.
\begin{definition}
\label{xi}
Let $({\cal A},{\cal W})$ and $({\cal A}',{\cal W}')$ be Schur pairs
 over $S$ (resp. $S$') of rank $r$ (resp. $r'$). Then a homomorphism
 $(\alpha,\xi): ({\cal A}',{\cal W}')\rightarrow ({\cal A},{\cal W})$ 
consists of
\begin{enumerate}
\item A morphism $\alpha: S \rightarrow S'$ such that 
$\alpha^{(*)}{\cal A}' \subseteq {\cal A}$;
\item A homomorphism 
$\xi \in {\cal H}om_{\Opower{S}{z}}(\Opower{S}{z}^{\oplus r'},
\Opower{S}{z}^{\oplus r})$
such that
for the induced homomorphism $\xi \in
{\cal H}om_{\Onegpower{S}{z}}(\Onegpower{S}{z}^{\oplus r'},
\Onegpower{S}{z}^{\oplus r})$ the inclusion $$\xi(\alpha^{(*)}{\cal W}') \subseteq 
{\cal W}$$
holds.
\end{enumerate}
\end{definition}
In this way, we get the category ${\frak S}$ of Schur pairs.
\begin{definition}
We define a full  subcategory ${\frak S}'$ of ${\frak S}$ as follows: 
$({\cal A},{\cal W}) 
\in \frak{S}'^{r}_{F}(S)$ if and only if $({\cal A},{\cal W}) 
\in\frak{S}^{r}_{F}(S)$ and ${\cal A}
\cap \Opower{S}{z} = {\cal O}_{S}$.
\end{definition}
The sense of this definition will become clear later on. For the time being, 
take it simply as a notation.
\section{Families of curves and sheaves}

In this section, we fix the geometric objects that we want to investigate,  
and we prove some basic 
properties of them. 
 The definitions we are going to make might seem a little technical. That is 
why much  space is given to 
 illustrations and examples. Even more examples may be found in 
 Chapter \ref{APPL}. 

Our first aim is to  study sheaves over families of curves. As to that, we 
need to fix three objects, 
namely: the base scheme, the total space and a sheaf on the total space.

\subsection{Families of curves}
As base schemes $S$ we allow all locally noetherian schemes.
\begin{definition}
\label{total C}
By a {\em pointed relative curve over $S$} we understand a scheme $C$ 
together with a locally projective
 morphism $\pi:C\rightarrow S$ and a section $P\subset C$ of $\pi$ such 
that the following holds:
\begin{enumerate}
\item $P$ is a relatively ample Cartier divisor in $C$.
\item For the sheaf ${\cal I}:= {\cal I}_{P}$ defining $P$ in $C$, 
${\cal I}/{\cal I}^{2}$ is a free 
${\cal O}_{P}$-module of rank 1.
\item Let $\widehat{\cal O}_{C}$ denote the formal completion of 
${\cal O}_{C}$ with respect to the 
ideal $\cal I$. Then $\widehat{\cal O}_{C}$ is isomorphic to 
$\Opower{P}{z}$ as a formal ${\cal O}_{P}$-algebra.
\item 
$\bigcap_{n\geq 0}\pi_{*}{\cal O}_{C}(-nP)=(0)$. 
\end{enumerate}
\end{definition}

Let us include here a couple of remarks and examples.
\freeline\\
{\bf Remark}
\begin{itemize}
\item Since $P$ is a section, $\pi|P : P \rightarrow S$ is an isomorphism. 
The sheaves 
${\cal I}/{\cal I}^{2}$ and $\widehat{\cal O}_{C}$ have their support in $P$. 
Consequently,
 Condition 2 is equivalent to $\pi_{*}({\cal I}/{\cal I}^{2})\cong {\cal O}_{S}$, 
while Condition 3 
  translates into: $\pi_{*}\widehat{\cal O}_{C}$ is isomorphic to $\Opower{S}{z}$ as 
a formal ${\cal O}_{S}$-algebra.

\item
Condition 4 is equivalent to the fact that, for every integer 
$n\in\Bbb{Z}$, the natural map
$$\pi_{*}{\cal O}_{C}(nP) \rightarrow \pi_{*}\widehat{{\cal O}_{C}(nP)}$$
is injective.
\end{itemize}
{\bf Example 1}\freespace
Let $S=Spec(k)$ for some field $k$. Then $C$ is a complete curve and $P$ 
corresponds to some smooth,
 $k$-rational point of $C$. The $k$-rationality is a consequence of the 
fact that $P$ is a section. Since
 $\widehat{\cal O}_{C,P}\cong \power{k}{z}$, the ring ${\cal O}_{C,P}$  
is regular. 

The  Condition 4 is satisfied if and only if the curve $C$ is reduced and 
irreducible.
\freeline\\
{\bf Example 2}\freespace
The motivation for studying curves over base schemes which are different 
from one point  comes mainly
 from the desire to investigate families of curves as considered in Example 1. 
Let us take, as an example,
  such a family over an integral $k$-scheme $S$. Since $S$ and the fibres of 
$\pi$ are irreducible, $C$ is 
  automatically irreducible. Let us assume, in addition, that $C$ is reduced. 
Then Condition 4 of the
   previous definition is satisfied. Conditions 2 and 3 amount to saying that 
our family is constant
    locally around the section $P$.
\freeline

Notice that the condition we impose on $P$ by assuming the triviality of  
${\cal I}/{\cal I}^{2}$ is very restrictive in the case where our base scheme 
$S$ is complete. This is expressed in the following proposition, which 
can be found in \cite{A1}:
\begin{proposition}
Let $f:X\rightarrow S$ be a family of nodal curves of genus $g > 0$ over 
a reduced and irreducible complete curve, and let $\Gamma \subset X$ be 
a section of $f$ not passing through any of the singular points of the 
fibres. Suppose that the general fibre of $f$ is smooth. Then
$$(\Gamma\cdot\Gamma)\leq 0.$$
Moreover, if $(\Gamma\cdot\Gamma)=0$, then the family $f:X\rightarrow S$, 
together with the section $\Gamma$, is an isotrivial family of 1-pointed 
nodal curves.
\end{proposition}

It is not hard to generalize the whole set-up to the case where 
${\cal I}/{\cal I}^{2}$ is locally free,
 but not free. One simply has to use $\prod_{n\geq 0}({\cal I}/{\cal I}^{2})^{n}$ 
instead of $\Opower{P}{z}$. 
 However, since we are mainly interested in local considerations, this generalization
 does not play such an
  important role that it would justify the technical effort.

Now, the first condition needs to be examined. The following lemma shows that it is 
almost automatically satisfied:
\begin{lemma}
\label{ample}
\begin{itemize}
\item If $S$ is irreducible, then $P$ is a Cartier divisor if and only if the conormal 
sheaf ${\cal I}/{\cal I}^{2}$ is a line bundle on $P$.
\item If  $P$ is a Cartier divisor on $C$, and if the morphism $\pi$ has irreducible fibres, 
then $P$  is relatively ample.
\end{itemize}
\end{lemma}
{\bf Remark}\freespace
The assumption on $P$ translates as follows: ${\cal I}={\cal I}_{P}$ is locally 
generated by one element 
which is not a zero divisor.
\freeline\\
{\bf Proof of the lemma}\freespace
The first statement is easy. Let us prove the second one. The question is local on $S$. So we
 are in the following situation:\\
$R$ is a noetherian ring, $\pi:C \rightarrow Spec(R)$ is a projective 
morphism, and $P\subset C$ is a section of $\pi$ and  an effective Cartier divisor. 

For all $n\in \Bbb{N}$, we have the following 
exact sequence:
$$ 0 \rightarrow {\cal O}_{C}((n-1)P) \rightarrow {\cal O}_{C}(nP) 
\rightarrow {\cal O}_{P}(nP) \rightarrow 0.$$
Since $P$ is affine, this induces a long exact sequence of cohomology groups:
$$
\begin{array}{ccccccccc}
0 & \rightarrow & H^{0}({\cal O}_{C}((n-1)P)) & \rightarrow & 
H^{0}({\cal O}_{C}(nP)) & 
\rightarrow & H^{0}({\cal O}_{P}(nP)) & {} & {}\\
{} & \rightarrow & H^{1}({\cal O}_{C}((n-1)P)) & \rightarrow & 
H^{1}({\cal O}_{C}(nP))   & \rightarrow & 0 .
\end{array}$$
Hence we have surjections $
H^{1}({\cal O}_{C}((n-1)P))) \surjection H^{1}({\cal O}_{C}(nP))$. Composing 
these surjections, we obtain, for each $n\in \Bbb{N}$, an epimorphism
$$\alpha_{n}:H^{1}({\cal O}_{C}) \surjection H^{1}({\cal O}_{C}(nP)).$$
Let $M_{n}$ denote the kernel of the mapping $\alpha_{n}$, which is an $R$ -  
submodule of $H^{1}({\cal O}_{C})$. From the very definition we get: $M_{n} 
\subseteq M_{n+1}$. But \cite{H1}, Thm. III.5.2., tells us that 
$H^{1}({\cal O}_{C})$ is a finitely generated $R$-module, hence noetherian. 
Therefore, there is an integer $N$ such that $M_{N} = M_{N+1}=\ldots$. This
implies that, for all $n > N$, the following sequence is  exact:
\begin{equation}
\label{noether}
0 \rightarrow H^{0}({\cal O}_{C}((n-1)P))) \rightarrow H^{0}({\cal O}_{C}(nP)) 
\rightarrow H^{0}({\cal O}_{P}(nP)) \rightarrow 0.
\end{equation}
But $P$ is affine. Therefore ${\cal O}_{P}(nP))$ is globally generated. This 
implies, using Nakayama's lemma and (\ref{noether}), that the global sections 
of ${\cal O}_{C}(nP)$  generate this sheaf in some neighborhood of $P$. 
By definition, base-points of the sheaf ${\cal O}_{C}(nP)$ are contained 
in $P$. Thus ${\cal O}_{C}(nP)$ itself has no base-points, i.e., it is globally 
generated. The sections of ${\cal O}_{C}(nP)$ define an $R$ - morphism 
$$\beta:C \rightarrow \Bbb{P}^{M}_{R}.$$
As $\beta$ is an $R$ - morphism, it is compatible with $\pi$, i.e.,  
$\pi= pr\circ \beta$, where $pr$ denote the natural projection from 
$\Bbb{P}^{M}_{R}$ onto $Spec(R)$. Since $P$ is a section, $\beta$ restricts 
to a closed embedding on $P$. In addition, we know that 
$\beta^{-1}(\beta(P)) = P$.  Now we want to prove that $\beta$ has finite 
fibres. Assume, on the  contrary, that there is a point $q\in \beta(C)$ such 
that $dim \beta^{-1}(q) \geq 1$. Let $X$ be an irreducible component of 
$\beta^{-1}(q)$ of dimension greater than 0. It is clear that $X$ does 
not intersect the divisor $P$. In particular, the intersection of $X$ 
with each fibre of $\pi$ is a closed subset of codimension at least 1. 
But the fibres of $\pi$ are assumed to be irreducible of dimension 1. 
Therefore $X$ cannot be contained in any fibre of $\pi$. But this is a 
contradiction, since $\pi(\beta^{-1}(q)) = pr(q)$. 
 So $\beta$ is a quasi-finite morphism. 

By \cite{H1}, Cor. II.4.8., $\beta$ is proper. From the Stein-factorization 
theorem (cf. \cite{G1}, Cor. 4.3.3.) one knows that in this case $\beta$ is 
also finite. 

By construction,  ${\cal O}_{C}(nP) = 
\beta^{*}({\cal O}_{\Bbb{P}^{M}_{R}}(1))$. ${\cal O}_{\Bbb{P}^{M}_{R}}(1)$ 
induces a very ample line bundle  on $\beta(C)$. Since $\beta$ is a finite 
morphism, this implies that 
 ${\cal O}_{C}(nP)$  is ample, too (cf. \cite{H2}, Prop.I.4.4.).
\QED
\freeline\\
{\bf Remark}\freespace
We have seen that a family of curves whose fibres satisfy the conditions of Definition 
\ref{total C}  almost matches the definition already. 

However, the converse is definitely not true. Namely, if $C$ is as in Definition 
\ref{total C}, 
then it may occur that some fibres of $\pi$ are not even integral
 (see Section \ref{base change} for examples).
\freeline\\
After this illustration we return to our general definition. 
\begin{lemma}
\label{affine covering}
Let $C$ be as in Definition \ref{total C}. Then:
\begin{enumerate}
\item $C$ is locally noetherian.
\item For each open affine subset $U$ of $S$, $\pi^{-1}(U)\setminus P$ is affine.
\item $S$ can be covered by open affine sets, 
$U_{i} = Spec(R_{i})$, such that for each $i$ there is an open affine 
subset $V_{i} = Spec(B_{i})$ of $\pi^{-1}(U_{i})$ containing $P\cap \pi^{-1}(U_{i})$.
\end{enumerate}
\end{lemma}
{\bf Proof}\freespace
\begin{enumerate}
\item $\pi$ is locally projective, hence it is locally of finite type. Therefore, 
together with $S$,  $C$ is also locally noetherian.
\item This is a consequence of the relative ampleness of $P$.
\item $\pi$ is locally of finite type, i.e., we can choose an open covering of $S$ by 
affine sets $U_{k}'= Spec (R_{k}')$ such that, for all $k$, 
$\pi:\pi^{-1}(U'_{k}) \rightarrow U'_{k}$ is proper, and  there are 
finitely many open affine sets $V_{k,l}' = Spec(B_{k,l}')$ satisfying:
\begin{itemize}
\item $\pi^{-1}(U_{k}') = \bigcup_{l} V_{k,l}'$,
\item $B_{k,l}'$ is a finitely generated $R_{k}'$ - algebra.
\end{itemize}
Let $U' := U_{k}'$ for some $k$. Choose a point $Q\in P\cap \pi^{-1}(U')$. $Q$ is 
contained in  one of the $V_{k,l}'$ which we denote by $V'$ for short. If 
$V'$ contains $P\cap \pi^{-1}(U')$, then we are done. Now let us assume that $V'$ does not 
contain $P\cap \pi^{-1}(U')$. $P\cap \pi^{-1}(U')$ is a closed subset of $\pi^{-1}(U')$. From 
the closedness of $\pi$ we conclude that 
$\pi((P\cap \pi^{-1}(U'))\setminus V')$ is closed in $U'$, and that $V'$ contains 
$P\cap \pi^{-1}(U'\setminus \pi(P\setminus V'))$.  The open set 
$U'\setminus \pi(P\setminus V')$ can be covered by open affine sets of the kind 
$U'(f) := Spec((R_{k})_{f})$, for certain elements $f\in R_{k}$, and we 
see that over each of the $U'(f)$'s the affine set $V'(f) = V_{k,l}'(f) := 
Spec ((B_{k,l}')_{f})$ has exactly the required property.
\QED
\end{enumerate}

Now let us illustrate Condition 3 of Definition \ref{total C}.
\begin{lemma} \label{powerseries}
Assume that ${\cal I}/{\cal I}^{2}$ is trivial. Then $\widehat{\cal O}_{C}$ is 
isomorphic to $\Opower{P}{z}$ 
if and only if each section of ${\cal I}/{\cal I}^{2}\cong {\cal O}_{P}$ lifts 
to a section of $\widehat{\cal O}_{C}$. 
In particular, one may interpret $z$ as the lift of a generating section of 
${\cal I}/{\cal I}^{2}$.
\end{lemma}
{\bf Remark}\freespace
Consequently, $\widehat{\cal O}_{C}$ is isomorphic to $\Opower{P}{z}$ if and only 
if, for any integer $n\in\Bbb{N}$, 
$n\geq 2$, the map
$$H^{0}({\cal I}/{\cal I}^{n}) \rightarrow H^{0}({\cal I}/{\cal I}^{2})$$
is surjective.
\freeline\\
{\bf Proof of the lemma}\freespace
First assume that $\widehat{\cal O}_{C} \cong \Opower{P}{z}$. Then 
${\cal I}/{\cal I}^{2} \cong {\cal O}_{P} \cdot z$ 
and of course
$$H^{0}(\Opower{P}{z}) = \power{H^{0}({\cal O}_{P})}{z} \surjection 
H^{0}({\cal O}_{P})\cdot z.$$

Now let $z$ be a global section of $\widehat{\cal O}_{C}$ such that 
$z(mod ({\cal I}^{2}))$ generates ${\cal I}/{\cal I}^{2}$. 
The sheaf $\widehat{\cal O}_{C}$ is defined to be the limit taken over the 
projective system ${\cal O}_{C}/{\cal I}^{n}$.  
We restrict our consideration to an open affine set $U = Spec(R) \subseteq S$
 such that there is an open subset $V = Spec(B)$ of $\pi^{-1}(U)$ which 
contains $p(U)$ and which is so that $I :={\cal I}(V)$ is free.  Then we 
construct the required isomorphism of sheaves locally on $U$. 

Choose an element $b\in I$ such that $I=bB$ and $[b] = [z] \in I/I^{2}$. 
Now take $f \in B/I^{m}$ for some $m$. Identifying $R$ with $\pi^{*}(R)  
\subset B$, there is a uniquely determined element $f_{1} \in R$ such that
 $f-f_{1} \in I/I^{m}$. By assumption, $I/I^{2}$ equals  $[b]\cdot R$, and thus there
 is a $f_{2} \in R$ such that $f-f_{1}-f_{2}b(mod I^{m}) \in I^{2}/I^{m}$. 
Since $I^{n}/I^{n+1}$ is generated by $b^{n}(mod I^{n+1})$, we can continue 
this process and get well-defined maps $$B/I^{m} \longrightarrow R \oplus
 R\cdot b \ldots \oplus R\cdot b^{m-1}.$$
These maps  give rise to a homomorphism of formal $R$ - algebras:
$$ \widehat{B} \longrightarrow \power{R}{b},$$
where $\widehat{B}$ denotes the completion of $B$ with respect to the ideal 
$I$. By \cite{Mat1}, Thm.8.12., this is an isomorphism. 

The fixed global section $z$ of  
$\widehat{\cal O}_{C}$ restricts on $U$ to an element of  
$\widehat{B}$. Since $[z], [b]\in I/I^{2}$ coincide, the  
homomorphism constructed above maps $z$ to an element $b\cdot (1+\alpha)$  
for some $\alpha\in 
\power{R}{b}\cdot b$. Remark that all those elements  $1+\alpha$ are 
invertible in $\power{R}{b}$ (cf. the appendix). Therefore the formal power series rings 
\power{R}{b} and \power{R}{z} are naturally isomorphic. So we finally 
get a well-defined isomorphism of $R$ - algebras:
$$\rho : \widehat{B} \longrightarrow \power{R}{b} 
\stackrel{\sim}{\longrightarrow} \power{R}{z}.$$

Here, $\rho$ does not depend on the choice of the local lift $b$ of 
$[z]$. Therefore, these locally defined isomorphisms glue together. It follows from  
Lemma \ref{affine covering} that $S$ has a covering by sets $U$ as 
considered above. Thus the construction gives a well-defined 
isomorphism of sheaves 
$$\rho : \widehat{\cal O}_{C}
\stackrel{\sim}{\longrightarrow} \Opower{P}{z}.$$
\QED
\freeline

\subsection{Families of sheaves}
\begin{definition}
\label{sheaf F}
As sheaves $\cal F$ on $C$ we admit all coherent sheaves such that
\begin{enumerate}
\item The formal completion $\widehat{\cal F}$ of $\cal F$ along $P$ is
 a free $\widehat{\cal O}_{C}$-module.
\item $\bigcap_{n\geq 0}\pi_{*}{\cal F}(-nP) =(0)$.
\end{enumerate}
\end{definition}
{\bf Remark 1}\freespace
Again, the second condition is equivalent to:
$$\pi_{*}{\cal F}(nP) \hookrightarrow \pi_{*}\widehat{{\cal F}(nP)}, \quad
\textrm{for all } n\in \Bbb{Z}.$$
{\bf Examples}
\begin{itemize}
\item If $C$ is an integral scheme, then the second condition is satisfied if 
and only if $\cal F$ is torsion free. 
\item If $C$ is a complete curve, then  the first condition is  also satisfied
 for torsion free sheaves, since these 
are free in smooth points. 
\item For general $S$ and $C$ as in Definition \ref{total C} and for any vector
 bundle $\cal F$ on $C$, Condition 2 
 of Definition \ref{sheaf F} is satisfied.
\end{itemize}
{\bf Remark 2}\freespace
The isomorphism 
$$\rho : \widehat{\cal O}_{C}\stackrel{\sim}{\longrightarrow}
 \Opower{P}{z}$$
makes $\widehat{\cal F}$ into an $\Opower{P}{z}$-module. If $\cal F$ 
satisfies Condition 1, then there is an isomorphism of 
$\widehat{\cal O}_{C}$-modules
$$\Phi : \widehat{\cal F}\stackrel{\sim}{\rightarrow} 
\widehat{\cal O}_{C}^{\oplus s}$$
or, equivalently, an isomorphism of $\Opower{P}{z}$-modules
$$\rho\circ\Phi : \widehat{\cal F} \stackrel{\sim}{\rightarrow}
\Opower{P}{z}^{\oplus s}.$$
A natural question is: 
What do the properties of $\cal F$ imply in general? In order to get a flavor 
of what is happening, we prove the following lemma.
\begin{lemma}
\label{loc. free}
Let $\cal F$ be a coherent sheaf of rank $r$, as in Definition \ref{sheaf F}. 
Then the following holds:
\begin{enumerate}
\item[(  i)] The rank of $\widehat{\cal F}$ as an  
$\widehat{\cal O}_{C}$ - module equals to the rank of $\cal F$. In particular,
 there is an open, dense subset of $C$ on
 which $\cal F$ has constant rank.
\item[( ii)] ${\cal F}|P$ is free of rank $r$.
\item[(iii)] $\cal F$ is locally free in a neighborhood of  $P$.
\end{enumerate}
\end{lemma}
{\bf Proof}\freespace
(i) is a consequence of (iii). 
Let us turn to the proof of  (ii). We  assume that 
$\widehat{\cal F}\cong \widehat{\cal O}_{C}^{\oplus s}$, for 
some $s\in \Bbb{N}$. Thus
$$ 
\begin{array}{rcl}
{\cal F}|P & = & {\cal F}\otimes_{{\cal O}_{C}} {\cal O}_{C}/{\cal I}\\
{} & = & \widehat{\cal F}/{\cal I}\widehat{\cal F}\\
{} & = & (\widehat{\cal O}_{C}/{\cal I}\widehat{\cal O}_{C})^{\oplus s}\\
{} & = & ({\cal O}_{C}/{\cal I})^{\oplus s}\\
{} & = & {\cal O}_{P}^{\oplus s},
\end{array}
$$
i.e., ${\cal F}|P$ is free of rank $s$. 

Now let $x$ be a point of $P$. By Nakayama's lemma we get a 
surjection 
$${\cal O}_{C,x}^{\oplus s} \surjection {\cal F}_{x}.$$
Denote by $K$ the kernel of this morphism. Taking the formal completion 
along $P$ we get an exact sequence
$$0\rightarrow \widehat{K} \rightarrow \widehat{{\cal O}}_{C,x}^{\oplus s} 
\rightarrow \widehat{{\cal F}}_{x}\rightarrow 0.$$
By \cite{Mat1}, Thm.2.4.,  $\widehat{K}$ vanishes. On the other hand, we know 
that $\widehat{{\cal O}}_{C,x}$ is faithfully flat over ${\cal O}_{C,x}$ 
(cf. \cite{Mat1}, Thm.8.14.). This implies that $K=0$ and we conclude:
${\cal O}_{C,x}^{\oplus s} \cong {\cal F}_{x}$. This completes the proof 
of the lemma.
\QED
\freeline\\
{\bf Remark}\freespace
In Section \ref{Eogd} we will see that, under certain conditions, the properties
(ii)
and (iii) of Lemma \ref{loc. free} are equivalent to Condition 1 in Definition 
\ref{sheaf F}. This will give us an
 explicit method to construct examples. 
A first step in this direction is the following lemma.
\begin{lemma}
\label{COMP.F}
Assume that $C$ is as in Definition \ref{total C}. Assume furthermore that 
$\cal F$ is locally free of rank $r$ in 
some neighborhood of  the section $P$ and that the restriction of $\cal F$ 
to $P$ is free. Then $\widehat{\cal F}$ is a
 free $\widehat{\cal O}_{C}$-module if and only if, for all $n\in\Bbb{N}$, 
the maps
$$H^{0}({\cal F}/{\cal I}^{n}\otimes {\cal F}) \rightarrow 
H^{0}({\cal F}/{\cal I}\otimes {\cal F})$$
are surjective.
\QED
\end{lemma}
\subsection{Definition of geometric data}
To end  this chapter, let us precisely define  the geometric objects which we 
want to relate to Schur pairs. Let $S$
 be a locally noetherian scheme.  
\begin{definition}
\label{def data}
By a {\em geometric datum of rank $r$ and index $F$ over $S$}, we mean 
a tupel
 $$(C,\pi,S,P,\rho,{\cal F},\Phi)$$ such that 
\begin{enumerate}
\item $C$ is a scheme.
\item $\pi:C\rightarrow S$ is a locally projective morphism.
\item $P \subset C$ is a section of $\pi$ such that 
\begin{itemize}
\item $P$ is a relatively ample Cartier divisor in $C$.
\item For the sheaf ${\cal I}:= {\cal I}_{P}$ defining $P$ in $C$, 
${\cal I}/{\cal I}^{2}$ is a trivial line bundle on $P$.
\item Let $\widehat{\cal O}_{C}$ denote the formal completion of 
${\cal O}_{S}$ with respect to the ideal $\cal I$. Then 
$\widehat{\cal O}_{C}$ is isomorphic to $\Opower{P}{z}$ as a formal 
${\cal O}_{P}$-algebra.
\item 
$\bigcap_{n\geq 0}\pi_{*}{\cal O}_{C}(-nP)=(0)$.
\end{itemize}
\item $\rho : \widehat{\cal O}_{C} 
\stackrel{\sim}{\longrightarrow}
 \Opower{P}{z}$ is an isomorphism of formal ${\cal O}_{P}$-algebras.
\item $\cal F$ is a coherent sheaf of rank $r$ on $C$ 
such that 
\begin{itemize}
\item The formal completion $\widehat{\cal F}$ of $\cal F$ along $P$ is a 
free $\widehat{\cal O}_{C}$-module of rank $r$.
\item $\bigcap_{n\geq 0}\pi_{*}{\cal F}(-nP) =(0)$.
\item $F = \gamma(\pi_{*}{\cal F}) - \gamma(R^{1}\pi_{*}{\cal F})\in K(S)$.
\end{itemize}
\item $\Phi : \widehat{\cal F}
\stackrel{\sim}{\longrightarrow} 
\widehat{\cal O}_{C}^{\oplus r}$ is an isomorphism of sheaves of 
$\widehat{\cal O}_{C}$ - modules.
\end{enumerate}
\end{definition}

\begin{definition}
\label{ident}
Two geometric data 
$$(C,\pi,S,P,\rho,{\cal F},\Phi) \textrm{ and } 
(C',\pi',S',P',\rho',{\cal F}',\Phi')$$ are identified if and only if 
\begin{itemize}
\item $S=S'$;
\item There is an isomorphism $\beta:C\rightarrow C'$ such that
\begin{itemize}
\item The diagram 
$$
\begin{array}{ccccc}
P  \subset& C & \stackrel{\beta}{\rightarrow} & C'&\supset  P'\\
\sim \searrow&\pi\downarrow&&\downarrow\pi'&\swarrow \sim\\
& S &=&S
\end{array}
$$
is commutative;
\item $\rho = \widehat{\beta}^{*}(\rho')$;
\end{itemize}
\item There is an isomorphism $\Psi: \beta^{*}{\cal F}' 
\rightarrow {\cal F}$ such that $\widehat{\beta}^{*}(\Phi') = \Phi\circ \widehat{\Psi}$.
\end{itemize}
In the sequel, we denote by ${\frak D}^{r}_{F}(S)$ the set of equivalence classes of 
geometric data of rank $r$ and index $F$ over $S$.
\end{definition}
\begin{definition}
\label{homo data}
A {\em homomorphism of geometric data} is a collection 
$$(\alpha,\beta,\Psi):(C,\pi,S,P,\rho,{\cal F},\Phi)\rightarrow
(C',\pi',S',P',\rho',{\cal F}',\Phi')$$
consisting of
\begin{enumerate}
\item A morphism $\alpha:S\rightarrow S'$;
\item A morphism $\beta:C\rightarrow C'$ such that
\begin{enumerate}
\item The following diagram is commutative:
$$
\begin{array}{ccccc}
P  \subset &C & \stackrel{\beta}{\rightarrow} & C'&\supset P' \\
\sim\searrow&\pi \downarrow &&\downarrow \pi'&\swarrow\sim\\
&S & \stackrel{\alpha}{\rightarrow} & S' 
\end{array}
$$
\item $\beta^{*}(P') = P$ as Cartier divisors; 
\item $\rho=\widehat{\beta}^{*}(\rho')$;
\end{enumerate}
\item A homomorphism of sheaves 
$\Psi: \beta^{*}{\cal F}' \rightarrow {\cal F}$.
\end{enumerate}

Two homomorphisms are identified iff they differ only by an identification 
isomorphism as defined in the previous definition.
\end{definition}
This establishes the category $\frak D$ of geometric data.
\freeline\\
Definition \ref{homo data} requires some justification.
\begin{lemma}
Let $(C,\pi,S,P,\rho,{\cal F},\Phi)$ and 
$(C',\pi',S',P',\rho',{\cal F}',\Phi')$ be geometric data and  
$\alpha:S\rightarrow S'$ and $\beta:C\rightarrow C'$
 morphisms such that Conditions 2(a) and 2(b) of Definition \ref{homo data} are satisfied. Then
\begin{itemize}
\item $\widehat{\beta}^{*}\widehat{\cal O}_{C'} \cong \widehat{\cal O}_{C}$ and 
\item $\widehat{\beta}^{*}{\widehat{{\cal F}}'} \cong \widehat{\beta^{*}{\cal
F}'}$, \end{itemize}
i.e., Condition 2(c) is well-formulated and for a homomorphism of sheaves
$\Psi : \beta^{*}{\cal F}' \rightarrow  {\cal F}$, the composition
$\pi_{*}((\rho\circ\Phi)\circ\widehat{\Psi}\circ(\widehat{\beta}^{*}
(\rho'\circ\Phi'))^{-1})$ belongs to 
${\cal H}om_{\power{{\cal O}_{S}}{z}}(\power{{\cal O}_{S}}{z}^{\oplus r'}, 
\Opower{S}{z}^{\oplus r})$.
\end{lemma}
{\bf Proof}\freespace
We consider the exact sequence of ${\cal O}_{C'}$-modules
\begin{equation}
\label{just 1}
0\rightarrow {\cal I}'^{n}\rightarrow {\cal O}_{C'} \rightarrow 
{\cal O}_{C'}/{\cal I}'^{n} \rightarrow 0.
\end{equation}
The support of the sheaf ${\cal O}_{C'}/{\cal I}'^{n}$ is contained in 
$P$, and, as a sheaf on $P$, 
${\cal O}_{C'}/{\cal I}'^{n}$ is known to be free. Therefore we can conclude, 
using \cite{Mat1}, ch.18, Lemma 2, that
$${\cal T}or^{{\cal O}_{C'}}_{1}({\cal O}_{C'}/{\cal I}'^{n}, {\cal O}_{C}) 
\cong {\cal T}or^{{\cal O}_{P'}}_{1}({\cal O}_{C'}/{\cal I}'^{n}, {\cal O}_{P})
=(0).$$
Consequently, the sequence
$$
0\rightarrow \beta^{*}{\cal I}'^{n}\rightarrow \beta^{*}{\cal O}_{C'} \rightarrow 
\beta^{*}({\cal O}_{C'}/{\cal I}'^{n}) \rightarrow 0
$$
is exact, as well. Using this we start our calculation:
$$
\begin{array}{rcl}
\widehat{\beta}^{*}\widehat{\cal O}_{C'} & = & 
lim_{n\to\infty} \beta^{*}({\cal O}_{C'}/{\cal I}'^{n})\\
&\cong& lim_{n\to\infty} \beta^{*}{\cal O}_{C'}/\beta^{*}{\cal I}'^{n}\\
&\cong &lim_{n\to\infty}{\cal O}_{C}/{\cal I}^{n}\\
& = & \widehat{\cal O}_{C}.
\end{array}
$$
We proceed analogously for the sheaf ${\cal F}'$. The sequence (\ref{just 1}) 
stays exact after 
tensoring with ${\cal F}'$, since we assumed ${\cal F}'$ to be locally free near 
$P'$. With the 
same conclusions as above  we finally get an exact sequence
$$
0\rightarrow \beta^{*}({\cal F}'\otimes{\cal I}'^{n})\rightarrow \beta^{*}({\cal F}') 
\rightarrow 
\beta^{*}({\cal F}'\otimes({\cal O}_{C'}/{\cal I}'^{n})) \rightarrow 0
$$
and, of course, ${\cal F}'\otimes({\cal O}_{C'}/{\cal I}'^{n}) \cong {\cal F}'/({\cal F}'\otimes
 {\cal I}'^{n})$. Now we can calculate again
$$
\begin{array}{rcl}
\widehat{\beta}^{*}\widehat{{\cal F}'}& = & lim_{n\to\infty} 
\beta^{*}({\cal F}'/{\cal F}'\otimes {\cal I}'^{n})\\
& \cong & lim_{n\to\infty} \beta^{*}({\cal F}')/\beta^{*}({\cal F}'\otimes {\cal I}'^{n})\\
&\cong & lim_{n\to\infty} \beta^{*}({\cal F}')/(\beta^{*}({\cal F}')\otimes {\cal I}^{n})\\
& = & \widehat{\beta^{*}({\cal F}')}.
\end{array}
$$
\QED
\freeline

Let us include  one more definition.
\begin{definition}
We define a full subcategory ${\frak D}'$ of $\frak D$ as follows:\\ The objects of this 
category are the equivalence 
classes of geometric data $(C,\pi,S,P,\rho,{\cal F},\Phi)$ such that 
$$\pi_{*}{\cal O}_{C} = {\cal O}_{S}.$$
\end{definition}
This subcategory will play an important role in Chapter \ref{Fam. DO}.
\freeline\\
{\bf Remark}\freespace
\label{field2}
\begin{itemize}
\item
Assume that $S=Spec(k)$, for some field $k$. Then the geometric datum
 $(C,\pi,S,P,\rho,{\cal F},\Phi)$ reduces to  $(C,P,\rho,{\cal F},\Phi)$. 
Using the same method as in the remark on page \pageref{Grass Mulase}, we 
see that this datum corresponds to a {\em quintet}  defined by M.~Mulase 
 \cite{M1},  where 
$$\rho :  \widehat{\cal O}_{C} \hookrightarrow \power{k}{y}$$
decomposes into
$$\widehat{\cal O}_{C} \stackrel{\sim}{\rightarrow} 
\power{k}{y^{r}} \hookrightarrow \power{k}{y}.$$
\item
Now let $S$ be any  scheme, $(C,\pi,S,P,\rho,{\cal F},\Phi)$ a 
geometric datum and $s\in S$ a closed point. Then we can restrict everything to the 
fibre $C_{s}$ of $C$ over $s$ and get a collection 
$$(C_{s},\pi|C_{s},\{s\},P|\{s\},\rho|C_{s},{\cal F}|C_{s},\Phi|C_{s}).$$
Assume that this is also a geometric datum. Then the restriction defines  a morphism 
of geometric data
$$(C_{s},\pi|C_{s},\{s\},P|\{s\},\rho|C_{s},{\cal F}|C_{s},\Phi|C_{s})
\rightarrow 
(C,\pi,S,P,\rho,{\cal F},\Phi).$$
\item Another example of morphisms of geometric data is the following: Let 
$(C,\pi,S,P,\rho,{\cal F},\Phi)$ be a 
geometric datum and $\alpha:S'\rightarrow S$ a flat morphism. Then the fibre 
product gives rise to another geometric datum. 
This fact will be shown in Section \ref{base change}.
\end{itemize}
\section{The relative Krichever functor}

After having defined both sides we want to relate, let us start with the construction of 
 a bijective contravariant functor 
between the category of Schur pairs and the category of geometric data. 
Throughout the chapter let us assume that $S$ is a locally noetherian scheme.
\subsection{Constructing Schur pairs}
Assume we are given a geometric datum $(C,\pi,S,P,\rho,{\cal F},\Phi)$ of 
rank $r$ and index $F$. 
Let us start with  
\begin{lemma}
\label{Cartier}
For each integer $n$, the maps $\rho$ and $\Phi$ induce isomorphisms:
$$
\begin{array}{*{5}{c}}
\rho & : & \widehat{{\cal O}_{C}(n\cdot P)} & 
\stackrel{\sim}{\longrightarrow} & 
 \Opower{P}{z}\cdot z^{-n},\\
\Phi & : & \widehat{{\cal F}(n\cdot P)} & 
\stackrel{\sim}{\longrightarrow} & 
\widehat{{\cal O}_{C}(n\cdot P)}^{\oplus r}.
\end{array}
$$
\end{lemma}
{\bf Proof}\freespace
This is clear from the fact that
$${\cal I} \cong {\cal O}_{C}(-P).$$
\QED
\begin{lemma}
\label{twist up}
For any  open affine subset $U$ of $S$, the natural maps
\begin{equation}
\label{loc g}
\begin{array}{ccc}
H^{0}(\pi^{-1}(U), {\cal O}_{C}(n\cdot P)) & \rightarrow & 
H^{0}(\pi^{-1}(U)\setminus P, {\cal O}_{C})\\
H^{0}(\pi^{-1}(U), {\cal F}(n\cdot P)) & \rightarrow & 
H^{0}(\pi^{-1}(U)\setminus P, {\cal F})\\
\end{array}
\end{equation}
induce isomorphisms
$$
\begin{array}{cccr}
lim_{n\to\infty}H^{0}(\pi^{-1}(U), {\cal O}_{C}(n\cdot P)) & 
\stackrel{\sim}{\rightarrow} & H^{0}(\pi^{-1}(U)\setminus P, 
{\cal O}_{C})&{}\\
lim_{n\to\infty}H^{0}(\pi^{-1}(U), {\cal F}(n\cdot P)) & 
\stackrel{\sim}{\rightarrow} & H^{0}(\pi^{-1}(U)\setminus P, 
{\cal F})&.\\
\end{array}
$$
\end{lemma}
{\bf Proof}\freespace
The maps in (\ref{loc g}) are inclusions, because  $P$ is locally given by 
an element which is neither a zero 
divisor in ${\cal O}_{C}$ nor in $\cal F$. 
We know that $P$ is ample relative to $S$. 
Therefore, for sufficiently large $n\in \Bbb{N}$, 
${\cal O}_{\pi^{-1}(U)}(n\cdot P)$ has a nonconstant global section. 
Now the assertion is a direct consequence of \cite{H1}, Lemma II.5.14. 
\QED
\begin{definition}
For a given geometric datum  $(C,\pi,S,P,\rho,{\cal F},\Phi)$ of 
rank $r$ and index $F$, 
 we define 
$$({\cal A},{\cal W}) := \chi_{r,F}(C,\pi,S,P,\rho,{\cal F},\Phi)$$
as follows:
%
%
%Def of A and W
%
%
\begin{displaymath}
\begin{array}{rcl}
{\cal A}(U) & := & \pi_{*}(\rho) (H^{0}(\pi^{-1}(U)\setminus P, 
{\cal O}_{C}))\\
{} & {}= &\pi_{*}(\rho) (lim_{n\to\infty}  H^{0}(\pi^{-1}(U), 
{\cal O}_{C}(n\cdot P)))\\
{} & {}= & \pi_{*}(\rho)(lim_{n\to\infty}  H^{0}(U, 
\pi_{*}{\cal O}_{C}(n\cdot P)))\\
{} & {}\subset & \negpower{{\cal O}_{S}(U)}{z},\\
{} & {} & {}\\
{\cal W}(U) & := & \pi_{*}(\rho\circ\Phi)(H^{0}(\pi^{-1}(U)\setminus P, 
{\cal F}))\\
{} & {}= &\pi_{*}(\rho\circ\Phi) (lim_{n\to\infty} H^{0}(\pi^{-1}(U), 
{\cal F}(n\cdot P)))\\
{} & {}= & \pi_{*}(\rho\circ\Phi)(lim_{n\to\infty} H^{0}(U, 
\pi_{*}{\cal F}(n\cdot P)))\\
{} & {}\subset & \negpower{{\cal O}_{S}(U)}{z}^{\oplus r}.
\end{array}
\end{displaymath}
\end{definition}
%
%
%End of definition
%
{\bf Remark }\freespace
\label{field1}
Assume that $S=Spec(k)$ for some field $k$. Then a geometric datum
 $(C,\pi,S,P,\rho,{\cal F},$ $\Phi)$ gives us a quintet 
$(C,P,\rho,{\cal F},\Phi)$ as defined in \cite{M1}, and $({\cal A},{\cal W})$ 
is the corresponding Schur pair defined there.
\freeline

In analogy to the case of a curve over a field we now want to identify the
(relative) cohomology of ${\cal O}_{C}$ and $\cal F$ via the  
sheaves $\cal A$ and $\cal W$. 

At first, we will see that the pole order of a local section along $P$ 
is exactly the order of the corresponding formal power series:
%
%Lemma
%
%
\begin{proposition}
\label{order}
For all integers $n\in\Bbb{Z}$:
$$
\begin{array}{*{5}{c}}
{\cal A} & \cap & \power{{\cal O}_{S}}{z}\cdot z^{-n} & = & 
\pi_{*}(\rho)(\pi_{*}({\cal O}_{C}(n\cdot P))),\\
{\cal W} & \cap & \power{{\cal O}_{S}}{z}^{\oplus r}\cdot z^{-n} & = & 
\pi_{*}(\rho\circ\Phi)(\pi_{*}({\cal F}(n\cdot P))).
\end{array}
$$
\end{proposition}
{\bf Proof}\freespace
This is obviously a local property. Therefore let $S=Spec(R)$ be affine, $R$ a 
local ring, $P \subset V=Spec(B)\subset C$,
$B$ another local ring,   $I:={\cal I}_{P}(V)=b\cdot B$,  $b$ a non-zero divisor
in $B$, and $M:=H^{0}(V,{\cal F})$. The composition of the maps 
$$H^{0}(C, {\cal F}(nP)) \rightarrow H^{0}(V,{\cal F}(nP))\rightarrow 
H^{0}(\widehat{{\cal F}(nP)})$$
is assumed to be injective. Therefore the first map must be injective. The
 second one is injective, 
a priori, since $B$ is a local ring. So we only have to show that
$$M_{b} \cap \widehat{M} = M.$$
Of course, we are done if we can show:
$$M \cap b\widehat{M} = bM,$$
i.e., $$M\cap I\widehat{M} = IM.$$
But this clear: Take $m\in M$ and consider its image in $\widehat{M}$. It can be 
identified with the 
sequence $(m(mod I^{n}M))_{n\in\Bbb{N}}$. This sequence belongs to $I\widehat{M}$
 if and only if, 
for all $n$, there are elements $i_{n}\in IM$ such that $m-i_{n}\in I^{n}M$. This 
 implies $m\in IM$ immediately.

Of course, the statement concerning $\cal A$ is proved in  the same way.
\QED
\freeline

Now we want to prove: 
%
%
%Theorem
%
%
\begin{proposition}
\label{cohomology}
There are isomorphisms of sheaves of ${\cal O}_{S}$ - modules
$$
\begin{array}{*{4}{c}}
{}&\frac{\negpower{{\cal O}_{S}}{z}}{{\displaystyle {\cal A}} + 
\power{{\cal O}_{S}}{z}} & \cong & 
R^{1}\pi_{*}{\cal O}_{C}\\
{} & {} & {}\\
\textrm{and }&\frac{\negpower{{\cal O}_{S}}{z}^{\oplus r}}{{\displaystyle {\cal W}} 
+ \power{{\cal O}_{S}}{z}^{\oplus r}} & \cong & 
R^{1}\pi_{*}{\cal F}.
\end{array}
$$
\end{proposition}
For the proof we need the following lemma.
\begin{lemma}
\label{iso completion}
Assume  that 
$U=Spec(R)$ is an open affine subset of $S$, and $V=Spec(B) \subseteq 
\pi^{-1}(U)$ is an open affine set containing $P\cap \pi^{-1}(U)$  such that, on
$V$, the divisor $P$ is
 given by a single element $b\in B$ and the restriction of $\cal F$ to $V$ is free. 
Then the natural maps
$$
\begin{array}{ccc}
H^{0}(V,{\cal O}_{C}) &\rightarrow &H^{0}(V,\widehat{\cal O}_{C}),\\
H^{0}(V,{\cal F})& \rightarrow &H^{0}(V,\widehat{\cal F})
\end{array}
$$
induce isomorphisms of $R$ - modules:
{\Large
$$
\begin{array}{*{5}{c}}
&\frac{H^{0}(V\setminus P,{\cal O}_{C})}{H^{0}(V,{\cal O}_{C})} & 
\stackrel{\sim}{\longrightarrow} & 
\frac{H^{0}(V\setminus P,\widehat{\cal O}_{C})}{H^{0}
(V,\widehat{\cal O}_{C})}&\\
{\normalsize\textrm{and }}&\frac{H^{0}(V\setminus P,{\cal F})}{H^{0}(V,{\cal F})} & 
\stackrel{\sim}{\longrightarrow} & 
\frac{H^{0}(V\setminus P,\widehat{\cal F})}{H^{0}
(V,\widehat{\cal F})}&.
\end{array}
$$} 
\end{lemma}
{\bf Proof}\freespace
Let us start with the investigation of the structure sheaf. 
The map 
$$B\rightarrow \widehat{B}$$
must be injective, since all elements of $B$ are locally given by quotients of 
elements of \linebreak
$\bigoplus_{n\geq 0} H^{0}(\pi^{-1}(U), {\cal O}_{C}(nP))$. 
Let us fix a natural number $n$. 
 We obtain a diagram of embeddings:
\begin{equation}
\begin{array}{*{8}{c}}
B & = & H^{0}(V,{\cal O}_{C}) & \rightarrow & H^{0}(V,\widehat{\cal O}_{C}) 
& = & \widehat{B} &{}\\
{} & {} & \downarrow & {} & \downarrow&{}\\
\frac{1}{b^{n}}B & = & H^{0}(V,{\cal O}_{C}(n\cdot P)) & \rightarrow & 
H^{0}(V,\widehat{{\cal O}_{C}(n\cdot P))} & = & \frac{1}{b^{n}}
\widehat{B}&.
\end{array}
\end{equation} 
Now the proof of Proposition \ref{order} implies that 
$$
\widehat{B}\cap \frac{1}{b^{n}} B = B.
$$
{}From this we obtain a natural inclusion:
$$\frac{\frac{1}{b^{n}} B}{B}\hookrightarrow 
\frac{\frac{1}{b^{n}} \widehat{B}}{\widehat{B}}.
$$
Taking the limit over $n$, we get  
 a monomorphism of $R$ - modules:
$$
\frac{B_{b}}{B} = 
\frac{{\displaystyle H^{0}(V\setminus P,{\cal O}_{C})}}{{\displaystyle 
H^{0}(V,{\cal O}_{C})}} 
\hookrightarrow 
lim_{n\to \infty}\frac{{\displaystyle H^{0}(V,\widehat{{\cal O}_{C}
(n\cdot P)})}}{{\displaystyle H^{0}(V,\widehat{\cal O}_{C})}}
= \frac{(\widehat{B})_{b}}{\widehat{B}}.
$$
We claim that this is, in fact, an isomorphism. Of course, we have $\poly{R}{b}\subseteq B$. 
So we obtain a diagram of
 inclusions:
$$
\begin{array}{ccc}
\poly{R}{b} & \rightarrow & B\\
\downarrow&&\downarrow\\
\negpoly{R}{b} & \rightarrow & B_{b},
\end{array}
$$
where $\negpoly{R}{b}\cap B = \poly{R}{b}$, since $b$ is a non-zero divisor in $B$. 
Therefore we end up with a
 chain of inclusions:
$$
\frac{\negpoly{R}{b}}{\poly{R}{b}}
\hookrightarrow
\frac{{\displaystyle H^{0}(V\setminus P,{\cal O}_{C})}}{{\displaystyle 
H^{0}(V,{\cal O}_{C})}}\hookrightarrow
\frac{{\displaystyle H^{0}(V\setminus P,\widehat{\cal O}_{C})}}
{{\displaystyle H^{0}(V,\widehat{\cal O}_{C})}} = 
\frac{\negpower{R}{b}}{\power{R}{b}}.
$$
But the first and the last term are canonically isomorphic 
$R$ - modules. Therefore all the monomorphisms appearing 
above are really isomorphisms.

Now we turn our attention to $\cal F$. The isomorphism ${\cal F}|V \cong 
{\cal O}_{V}^{\oplus r}$ extends to
 an isomorphism $\widehat{\cal F}\cong \widehat{\cal O}_{V}^{\oplus r}$, which 
also respects the localization by
  $b$. So we get the claim concerning $\cal F$ by applying the according statement 
on ${\cal O}_{C}$.
\QED
\freeline

Let us proceed to the
\freeline\\
{\bf Proof of the Proposition \ref{cohomology}}\freespace
As the isomorphism constructed in Lemma \ref{iso completion} is natural, 
it is compatible with intersections of affine sets. 
Now,  once again, it is sufficient to consider an open affine set $U$ of $S$, as in 
Lemma \ref{iso completion}. 
Then $\pi^{-1}(U)$ is covered by the 
affine sets $V$ and $\pi^{-1}(U)\setminus P$. Since $\pi|\pi^{-1}(U)$ is 
a separated morphism over an affine scheme, we can apply \cite{H1}, Thm. 
III.4.5., and get, with the aid of  Lemmas \ref{twist up}
 and  \ref{iso completion}, 
$$
\begin{array}{rcl}
R^{1}\pi_{*}{\cal O}_{C}(U) & = & H^{1}(\pi^{-1}(U), {\cal O}_{C})\\
{}&{}&{}\\
{} & = & \frac{{\displaystyle H^{0}(V \setminus P, {\cal O}_{C})}}{{
\displaystyle H^{0}(V, {\cal O}_{C}) + H^{0}(\pi^{-1}(U)\setminus P,
 {\cal O}_{C})}}\\
{}&{}&{}\\
{} & = & \frac{{\displaystyle H^{0}(V\setminus P, 
\widehat{\cal O}_{C})}}{{\displaystyle H^{0}(V, \widehat{\cal O}_{C}) 
+ H^{0}(\pi^{-1}(U)\setminus P, {\cal O}_{C})}}\\
{}&{}&{}\\
{} & \cong & \frac{{\displaystyle \Onegpower{S}{z}(U)}}{{\displaystyle 
\Opower{S}{z}(U) + {\cal A}(U)}}.
\end{array}
$$
The proof of the statement concerning $\cal F$ can be given in  the same way.
\QED
\begin{corollary}
If $(C,\pi,S,P,\rho,{\cal F},\Phi)$ is a geometric datum of rank $r$ and 
index $F$, then 
$({\cal A},{\cal W}) = \chi_{r,F}(C,\pi,S,P,\rho,{\cal F},\Phi)$
is a Schur pair of rank $r$ and index $F$ over $S$.
\end{corollary}
{\bf Proof}\freespace
By construction, $\cal A$ is a quasicoherent sheaf of ${\cal O}_{S}$ - 
subalgebras of $\Onegpower{S}{z}$, $\cal W$ is a quasicoherent sheaf of 
${\cal O}_{S}$ - modules, and ${\cal A}\cdot{\cal W}\subseteq{\cal W}$. 
The fact that $\cal A$ and $\cal W$ are elements of the infinite
 Grassmannians of rank 1 and $r$, respectively, follows from  
Propositions \ref{order} and 
\ref{cohomology}, and the fact that  $R^{i}\pi_{*}{\cal O}_{C}$ and $R^{i}\pi_{*}{\cal F}$ 
are coherent sheaves for locally projective 
morphisms $\pi$ and  for all $i$ (cf. \cite{H1}, Thm.III.8.8).
\QED
\begin{corollary}
\label{ident. Schur}
If $(C,\pi,S,P,\rho,{\cal F},\Phi)$ and 
$(C',\pi',S',P',\rho',{\cal F}',\Phi')$ 
are identified as geometric data, then we obtain identical
 corresponding Schur 
pairs $({\cal A}, {\cal W})$ and $({\cal A}', {\cal W}')$.
\end{corollary}
{\bf Proof}\freespace
This follows from an easy calculation. We use the isomorphisms $\beta$ 
and $\Psi$ and get:
$$
\begin{array}{rcl}
{\cal A} & := & \pi_{*}(\rho)(lim_{n\to\infty} 
\pi_{*}{\cal O}_{C}(n\cdot P))\\
& {}= & \pi_{*}(\widehat{\beta}^{*}(\rho'))(lim_{n\to\infty}  
\pi_{*}(\beta^{*}{\cal O}_{C'}(n\cdot P')))\\
& {}= & \pi'_{*}(\rho')(lim_{n\to\infty} 
\pi'_{*}{\cal O}_{C'}(n\cdot P'))\\
& {}= & {\cal A}'.
\end{array}
$$
Analogously:
$$
\begin{array}{rcl}
{\cal W}' & := & \pi'_{*}(\rho'\circ \Phi')(lim_{n\to \infty}\pi'_{*}
{\cal F}'(n\cdot P'))\\
& {}= & \pi'_{*}(\rho'\circ \Phi')(lim_{n\to \infty} \pi_{*}\beta^{*}
{\cal F}'(n\cdot P'))\\
& {}= & \pi_{*}(\widehat{\beta}^{*}(\rho'\circ \Phi'))(lim_{n\to \infty} \pi_{*}((\beta^{*}
{\cal F}')(n\cdot P)))\\
& {}= & \pi_{*}(\rho\circ \Phi\circ\widehat{\Psi})(lim_{n\to \infty} \pi_{*}((\beta^{*}
{\cal F}')(n\cdot P)))\\
& {}= & \pi_{*}(\rho\circ \Phi)(lim_{n\to \infty} \pi_{*}(
{\cal F}(n\cdot P)))\\
& {}=&{\cal W}.
\end{array}
$$
\QED
\begin{proposition}
\label{morph}
Homomorphisms of geometric data induce homomorphisms of the 
corresponding Schur pairs.
\end{proposition}
{\bf Proof}\freespace
Let $(\alpha,\beta,\Psi):(C,\pi,S,P,\rho,{\cal F},\Phi)
\rightarrow
(C',\pi',S',P',\rho',{\cal F}',\Phi')$ be a homomorphism of 
geometric data. Let $({\cal A}, {\cal W})$ and $({\cal A}', 
{\cal W}')$ be the Schur pairs associated to the given geometric 
data. We want to construct a homomorphism of Schur pairs
$$(\alpha,\xi):({\cal A}', {\cal W}')\rightarrow ({\cal A}, 
{\cal W}).$$
Of course, the morphism $\alpha:S\rightarrow S'$ is taken 
directly from $(\alpha,\beta,\Psi)$, whereas $\xi$ is defined  as
$\pi_{*}((\rho\circ\Phi)\circ\widehat{\Psi}\circ(\widehat{\beta}^{*}
(\rho'\circ\Phi'))^{-1})$. Now we
apply the properties of the given morphism $\beta:C\rightarrow C'$ 
and derive:
$$
\begin{array}{rcl}
\alpha^{(*)}{\cal A}' & = & \alpha^{(*)}(\pi'_{*}(\rho')(lim_{n\to\infty}   
\pi'_{*}{\cal O}_{C'}(n\cdot P')))\\
{} & \subseteq & \pi_{*}(\widehat{\beta}^{*}(\rho'))(lim_{n\to\infty} \pi_{*}
\beta^{*}{\cal O}_{C'}(n\cdot P'))\\
{} & = & \pi_{*}(\rho)( lim_{n\to\infty} \pi_{*} 
{\cal O}_{C}(n\cdot P))\\
{} & = & {\cal A}.
\end{array}
$$
As for ${\cal W}$ and ${\cal W}'$, we obtain:
$$
\begin{array}{rcl}
\xi(\alpha^{(*)}{\cal W}') & = & \xi(\alpha^{(*)}(\pi'_{*}(\rho'\circ\Phi')(lim_{n\to\infty} 
\pi'_{*}{\cal F}'(n\cdot P'))))\\
&\subseteq& (\xi\circ\pi_{*}(\widehat{\beta}^{*}(\rho'\circ\Phi')))(lim_{n\to\infty} 
\pi_{*}\beta^{*}{\cal F}'(n\cdot P'))\\

&=& \pi_{*}(\rho\circ\Phi\circ\widehat{\Psi})(lim_{n\to\infty} 
\pi_{*}\beta^{*}{\cal F}'(n\cdot P'))\\
&\subseteq&\pi_{*}(\rho\circ\Phi)(lim_{n\to\infty} \pi_{*}{\cal F}(n\cdot P'))\\
&=& {\cal W}.
\end{array}
$$
So we really have constructed a homomorphism of Schur pairs.
\QED
\begin{corollary}
$\chi$ is a contravariant functor from the category of geometric 
data to the category of Schur pairs.
\QED
\end{corollary}
\begin{definition}
The functor $\chi$ is called the {\em Krichever functor}.
\end{definition}
\subsection{Constructing geometric data}
Now assume that we are given a Schur pair $({\cal A},{\cal W})$ of 
rank $r$ and index $F$ over the scheme $S$.  We start with some 
general observations.
\begin{lemma}
\label{A0}
\begin{enumerate}
\item For all $n\in\Bbb{Z}$, ${\cal A}^{(n)}$ and ${\cal W}^{(n)}$ are
  coherent sheaves. In particular, 
${\cal A}^{(0)}$ is a coherent sheaf of ${\cal O}_{S}$-algebras. 
\item If $U = Spec(R)$ is an open affine subset of $S$, then there is an 
integer $M\in\Bbb{Z}$ 
(possibly depending on $U$) such that
$${\cal A}^{(-M)}(U) = {\cal W}^{(-M)}(U) = (0).$$
\item All local sections of ${\cal A}^{(-1)}$ are nilpotent.
\end{enumerate}
\end{lemma}
{\bf Proof}\freespace
\begin{enumerate}
\item Let $U=Spec(R)$ be an affine open subset of $S$. We have to show that
 ${\cal A}^{(n)}(U)$ and 
${\cal W}^{(n)}(U)$ are  finitely generated. ${\cal A}^{(0)}(U)$ and 
${\cal W}^{(0)}(U)$ are finitely
 generated by the definition of Schur pairs. Since $R$ is a noetherian 
ring, this immediately proves 
 the statement for all $n\leq 0$. Now assume that $n>0$. Then 
${\cal A}^{(n)}(U)/{\cal A}^{(0)}(U)$ is isomorphic to a submodule of 
$\power{R}{z}\cdot z^{-n}/\power{R}{z}$, hence finitely generated. But 
this already implies that
 ${\cal A}^{(n)}(U)$ itself is finitely generated. The same method of
 proof may be used for $\cal W$. 
\item Without loss of generality we prove the statement for $\cal A$. 
${\cal A}^{(0)}$ is generated by
 finitely many elements
$f_{1},\ldots,f_{m}$. Let us write
$$f_{i} = \sum_{j\geq 0} \lambda_{i,j} z^{j}.$$
The statement ${\cal A}^{(-M)}(U)\neq (0)$ is equivalent to:
 There are elements $\mu_{1},\ldots, \mu_{m} \in R$ such that 
$\sum_{i=1}^{m} \mu_{i}\lambda_{i,j} = 0$ for $j=0,\ldots, M-1$, but 
$\sum_{i=1}^{m} \mu_{i} f_{i} \neq 0$.

We denote by ${\cal N}_{l}$ the submodule of $R^{m}$ generated by the vectors 
$(\lambda_{\cdot ,0}), \ldots, (\lambda_{\cdot ,l})$. These modules form an 
ascending chain of submodules of $R^{m}$:
$${\cal N}_{0}\subseteq {\cal N}_{1}\subseteq{\cal N}_{2}\subseteq \ldots
\subseteq R^{m}.$$
Since $R^m$ is a noetherian module, there is an integer $T$ such that 
${\cal N}_{l}={\cal N}_{T}$, for all $l\geq T$. 

We claim that ${\cal A}^{(-T-1)}(U)=(0)$. If this were false, we could 
find elements \linebreak[4] 
$\mu_{1},\ldots, \mu_{m} \in R$ such that 
$\sum_{i=1}^{m} \mu_{i}\lambda_{i,j} = 0$, for $j=0,\ldots, T$. But by the 
definition of ${\cal N}_{l}$
 this already implies $\sum_{i=1}^{m} \mu_{i} f_{i} = 0$. So we are done. 
\item 
The last part of the lemma is easy. Assume that for some $U$, ${\cal A}^{(-1)}(U)$ 
would contain an element 
$f$ which is not nilpotent. Then, for all $n\in\Bbb{N}$,
$$0\neq f^{n} \in {\cal A}^{(-n)}(U),$$
and this is a contradiction.
\end{enumerate}
\QED
%
%
%Lemma
%
%
\begin{lemma}
\label{a and b}
Let $U = Spec(R)$ be an open affine subset of $S$. Then ${\cal A}(U)$ 
is a finitely generated $R$ - algebra of relative dimension 1.
\end{lemma}
{\bf Proof}\freespace
By assumption, 
$\frac{\negpower{R}{z}}{\power{R}{z} + {\displaystyle {\cal A}(U)}}$ is a 
finitely generated $R$ - module. So we can choose finitely many elements 
$b_{1},\ldots,b_{m}\in \negpower{R}{z}$ such that $[b_{1}],\ldots,[b_{m}]$ 
are generators of this module. Denote by $N$ the maximum of the orders of 
the $b_{j}$'s. We may assume that $N\geq 2$. Now it is straightforward that:
$$\negpower{R}{z} = \power{R}{z}\cdot z^{-N} + {\cal A}(U).$$
Therefore, for all $n>N$, there is a monic element of order $n$ in 
${\cal A}(U)$. Let us choose 
 monic elements $a, b\in {\cal A}(U)$ of order 
$2N+2$ and 
$2N+1$, respectively. We claim that
\begin{equation}
\label{a and b,2}
\negpower{R}{z} = \power{R}{z}\cdot z^{-(2N+1)(2N+2)} + \poly{R}{a,b}.
\end{equation}
In order  to prove this, we choose an integer $n \geq (2N+1)(2N+2)$. Then we can 
find $m \geq 0$ and $0 \leq l < 2N+2$ such that 
$$
\begin{array}{*{4}{c}}
n & = & (2N+1)(2N+2) + m\cdot (2N+2) + l &{}\\
{} & = & (2N+1)(2N+2) + m\cdot (2N+2) + l\cdot((2N+2) - (2N+1)) &{}\\
{} & = & (m + l)\cdot (2N+2) + (2N+2 - l)\cdot (2N+1) &{}\\
{} & = & ord( a^{m+l}\cdot b^{2N+2-l} ) &.
\end{array}
$$
Since $a$ and $b$ are monic, this proves (\ref{a and b,2}) and, in particular, 
the identity
\begin{equation}
{\cal A}(U) = {\cal A}(U)^{((2N+1)(2N+2))} + \poly{R}{a,b}.
\end{equation}
${\cal A}(U)^{((2N+1)(2N+2))}$ is a finitely generated $R$ - module 
(cf. Lemma \ref{A0}). This 
implies that  ${\cal A}(U)$ is a finitely generated $R$ - algebra of 
relative dimension at most 2. From the choice of $a$ and $b$, it is also 
clear that the relative dimension is greater than zero. 

Now we only need to prove that $a$ and $b$ satisfy a polynomial relation 
with coefficients in $R$. Obviously, $b$ does not lie in $\poly{R}{a}$. 
Let $\{u_{1},\ldots,u_{q}\}$ be a set of generators of 
$\poly{R}{a,b}^{((2N+1)(2N+2)-1)}$ as an $R$ - module, and put
$$v_{n} := a^{m+l}\cdot b^{2N+2-l}$$
for $n\geq (N+1)(N+2)$ and $n  = (2N+1)(2N+2) + m\cdot (2N+2) + l$ (see above). 
Then the set $\{u_{1},\ldots,u_{q}\}\cup \{v_{n}\}_{n\geq (2N+1)(2N+2)}$ generates the 
$R$ - module $\poly{R}{a,b}$. For all $M$, $a^{M}$ can be written as a linear 
combination of the $u_{j}$'s and $v_{j}$'s. Since none of the $v_{j}$'s is a 
power of $a$, and since there are only finitely many $u_{j}$'s,  the representing 
linear combination for a sufficiently 
high power of $a$ is exactly the required 
polynomial relation between $a$ and $b$. This completes the proof.
\QED
\begin{lemma}
\label{m+n}
Let $U=Spec(R)$ be an open affine subset of $S$ and let $N\in\Bbb{N}$ be 
such that $\negpower{R}{z} = 
\power{R}{z}\cdot z^{-N} + {\cal A}(U)$ (cf. Lemma \ref{a and b}). Then, 
for all $n,m \geq 2N+1$:
$${\cal A}(U)^{(m)}\cdot{\cal A}(U)^{(n)}={\cal A}(U)^{(m+n)}.$$
\end{lemma}
{\bf Proof}\freespace
Of course, ${\cal A}(U)^{(m)}\cdot{\cal A}(U)^{(n)}\subseteq {\cal A}(U)^{(m+n)}$. 
To see the other
 inclusion it is sufficient to show that ${\cal A}(U)^{(m)}\cdot{\cal A}(U)^{(n)}$ 
contains monic elements of the orders 
$m+1, \ldots, m+n$. 

First take $1\leq i\leq n-N$. Then we can split $m+i = (m-N) + (N+i)$. Since $m-N$ is 
greater than $N$,
 by assumption, ${\cal A}(U)^{(m)}$ contains a monic element of order $m-N$. On
 the other hand, 
 $N+1\leq N+i \leq n$. Therefore, ${\cal A}(U)^{(n)}$ contains a monic element 
of order $N+i$. 

As a second case consider now the terms $m+i$ for $n-N<i\leq n$. In this case, 
$i$ is greater than 
$N$, so ${\cal A}(U)^{(n)}$ contains a monic element of order $i$, and, of course, 
${\cal A}(U)^{(m)}$ 
contains a monic element of order $m$. This completes the proof.
\QED
\freeline

Now we aim at defining geometric objects from a given Schur pair. First, we define 
a sheaf  of graded ${\cal O}_{S}$ - algebras as follows:
$$grd({\cal A}) := {\cal O}_{S}\oplus\bigoplus_{n\geq 1} {\cal A}^{(n)},$$
i.e., $grd({\cal A})_{0} = {\cal O}_{S}$ and  $grd({\cal A})_{n} = 
{\cal A}^{(n)}$, for $n\geq 1$. 
Now, we define a scheme $C$ by
$$C := Proj(grd({\cal A})).$$
This scheme comes equipped with a projection morphism $\pi$ to $S$, and 
Lemma \ref{a and b} says that $C$ is a curve over $S$. 
\freeline\\
{\bf Remark}\freespace
\begin{itemize}
\item The morphism $\pi$ factors over the scheme $Spec({\cal A}^{(0)})$. 
By Lemma \ref{A0}, this scheme 
 is finite over $S$, and $Spec({\cal A}^{(0)})_{red}=S_{red}$.
\item
It is well-known that, for any $m\in \Bbb{N}$, the scheme $C$ is naturally 
isomorphic to
$$C^{(m)} := Proj({\cal O}_{S}\oplus\bigoplus_{n\geq 1} {\cal A}^{(mn)}).$$
On the other hand, Lemma \ref{m+n} implies that for any affine open subset $U\subseteq S$ 
 there is a number $m$ such that 
${\cal O}_{S}(U)\oplus\bigoplus_{n\geq 1} {\cal A}^{(mn)}(U)$ is generated by 
${\cal A}^{(m)}(U)$ 
as an ${\cal O}_{S}(U)$ - algebra. Moreover, we know that, for
all $m$,  ${\cal A}^{(m)}$ is coherent (see  Lemma \ref{A0}). 
Then we obtain from the general theory developed in \cite{H1}, II.7., that
 $C$ is locally projective over $S$. 
\item The localization of $grd({\cal A})$ by the section $1\in grd({\cal A})_{1}$ 
can be identified with
 $\cal A$, thus $C$ contains the (relatively) affine subset $Spec({\cal A})$.
\end{itemize}
\freeline 

Analogously, we define
$$grd({\cal W}) := \bigoplus_{n\in \Bbb{Z}} {\cal W}^{(n)}.$$

This gives us a locally finitely generated $grd({\cal A})$ - module. The sheaf 
$\cal F$ on $C$ is defined to be
$${\cal F} := (grd({\cal W}))^{\sim}.$$

We would like to see that the curve $C$ just constructed is exactly a curve of the 
type we started with.
\begin{theorem}
\label{Schur curve}
There is a section $P\subset C$ of $\pi$ 
such that $C\setminus P$ is precisely $Spec({\cal A})$.
 $P$ is a relatively ample Cartier divisor, its conormal 
sheaf  is free of rank 1 on $P$, and  $\widehat{\cal O}_{C}$ is isomorphic to 
$\Opower{P}{z}$. Finally, 
$\bigcap_{n\geq 0} \pi_{*}{\cal O}_{C}(-nP) = (0)$. 

\end{theorem}
{\bf Proof}\freespace
Again, we restrict everything to an affine open subset $U= Spec(R)$ of $S$. Let 
$a$ and $b$ be the elements we constructed in Lemma \ref{a and b}. We can 
view $a$ as an element of $(grd({\cal A}))(U)_{2N+2}$. Let us localize 
$grd({\cal A})$ by $a$:
$$
\begin{array}{*{5}{c}}
B & := & (grd({\cal A}))(U)_{(a)} & = & 
\left\{ \frac{g}{a^{n}}/ n\in \Bbb{N}, g \in (grd({\cal A}))(U)_{(2N+2)\cdot n} 
\right\}\\
&&{} & = & \left\{ \frac{g}{a^{n}}/ n\in \Bbb{N}, g \in 
{\cal A}(U)^{((2N+2)\cdot n)} \right\}\\
&&{} & \subseteq & \power{R}{z}.
\end{array}
$$
Especially, $y := \frac{b}{a}$ gives us a monic element of 
$B$ of order -1. But of course, for this 
element $y$ we obtain:
$$\power{R}{y} \cong \power{R}{z},$$
since all elements of order zero with invertible leading term are 
invertible in \power{R}{z}. Thus, we have to consider the situation:
\begin{equation}
\label{ideal in local}
\poly{R}{y} \subseteq B \subseteq \power{R}{y} = 
\power{R}{z} . 
\end{equation}
Let $I :=  B \cap \power{R}{z}\cdot z$. Obviously, 
this is an  ideal of $B$. Let $P$ be the closed subscheme of $Spec(B)$ defined
 by $I$. One easily sees from 
(\ref{ideal in local}) that 
$B/I$ is naturally isomorphic to $R$. Therefore, $P$ is a section of the 
projection morphism $\pi$.  Observe 
that the definition of $P$ does not depend on the choice of $a$ 
and $b$. For the ideal sheaf $\cal I$ of $P$, we get immediately:
$${\cal I}/{\cal I}^{2} = [y] \cdot ({\cal O}_{C}/{\cal I}) \cong {\cal O}_{P}\cdot z,$$
i.e., ${\cal I}/{\cal I}^{2}$ is free of rank 1.

Now let us prove that $C\setminus P= Spec({\cal A})$. We do this again on the 
affine open subset 
$U=Spec(R)$ of $S$. In a first step, we restrict our consideration to the open 
set $D_{+}(a)=Spec(B)$. 
$P$ is contained in this set,  and $D_{+}(a)$ and $Spec({\cal A}(U))$ cover
 $\pi^{-1}(U)$. Assume, on the
  contrary, that there is a graded prime ideal ${\frak p} \subset grd({\cal A}(U))$ 
such that
   $1\in {\frak p}_{1}$ and $a\in {\frak p}_{2N+2}$. By the choice of $a$, this 
already implies that 
   ${\frak p}_{n} = grd({\cal A}(U))_{n}$ for all sufficiently large $n$ 
(cf. Lemma \ref{m+n}). This is a contradiction.

On $D_{+}(a)= Spec(B)$, $Spec({\cal A}(U))$ is given as $D(\frac{1}{a})$. 
Now the second step is the following:
 Let ${\frak p}$ be a prime ideal of $B$ containing $\frac{1}{a}$. We want 
to show that $\frak p$ belongs to $P$.
  This holds if and only if $\frak p$ contains the ideal
$$\left\{ \frac{g}{a^{n}}/ n\in\Bbb{N}, ord(g)<ord(a)n \right\}.$$
Assume that $ord(g)<ord(a)n$ and set $m:=n\cdot ord(a)-ord(g)$.  Then:
$$\frac{g^{ord(a)}}{a^{n\cdot ord(a)}}=\frac{g^{ord(a)}}{a^{ord(g)+m}}
= \frac{g^{ord(a)}}{a^{ord(g)}}\cdot\frac{1}{a^{m}}\in {\frak p}.$$
But $\frak p$ was assumed to be prime, hence $\frac{g}{a^{n}}\in {\frak p}$. 
So we see that $\frac{1}{a}\cdot B$ defines the same closed subset of 
$D_{+}(a)$ as $I$.

As for  the last step, we have to prove that $D(\frac{1}{a})$ does not
 intersect $P$. Take a prime ideal
 ${\frak p}\in P$. It is sufficient to show that 
$$B_{1/a}\cdot {\frak p} = B_{1/a}= {\cal A}_{a}.$$
But this is obvious.

Our next aim is to show that the completion of ${\cal O}_{C}$ along $P$ is 
isomorphic to $\Opower{P}{z}$. 
To do this, we use the inclusions (\ref{ideal in local}). By \cite{Mat1},
 Thm. 8.1., we only need to prove
 that the $(y)$-adic topology on $\power{R}{y}$ induces the $I$-adic
 topology on $B$, and this one induces 
 the $(y)$-adic topology on $\poly{R}{y}$. Note that once we have shown 
the first fact, the second one follows immediately. 

We saw that $I$ and $\frac{1}{a}B$ define the same closed subset of $Spec(B)$. 
Therefore, both ideals define
 the same topology on $B$.
Since, obviously, 
$$(\frac{1}{a}B)^{n} \subseteq \power{R}{y}\cdot y^{n(2N+2)}\cap B,$$
there only remains to show that for each $k\in \Bbb{N}$ there is an integer 
$N(k)$ such that
$$\power{R}{y}\cdot y^{N(k)} \cap B \subseteq (\frac{1}{a}B)^{k}.$$
We claim that this is true for $N(k) = k\cdot (2N+2) = k\cdot (ord(a))$. 

{}From the definition of $B$ we get:
$$\power{R}{y}\cdot y^{k\cdot (ord(a))} \cap B = 
\{ \frac{g}{a^{\alpha}} / g\in {\cal A}(U), ord(g) \leq ord(a)(\alpha -k) \}.$$
Let us consider such an element $\frac{g}{a^{\alpha}}\in 
\power{R}{y}\cdot y^{k\cdot (ord(a))} \cap B$. 
The inequality $ord(g) \leq ord(a)(\alpha -k)$ implies 
$ord(g\cdot a^{k}) \leq ord(a)\alpha $. Therefore,
 $\frac{g\cdot a^{k}}{a^{\alpha}}$ is an element of $B$, 
i.e., $\frac{g}{a^{\alpha}} \in (\frac{1}{a}B)^{k}$.  

So we have shown that the $I$-adic completion of $B$ is isomorphic 
to $\power{R}{y} = \power{R}{z}$. This 
isomorphism obviously does not depend on the choice of $a$ and $b$. 
Therefore, it extends to an isomorphism
$$\rho: \widehat{\cal O}_{C} \stackrel{\sim}{\longrightarrow}
\Opower{P}{z}.$$

{}Furthermore, it is now an easy consequence of \cite{Mat1}, Thm. 7.5., 
that the ideal $I$ is locally free of
 rank 1. In fact, along $P$, this ideal is generated by the element $y$. 
This implies that $P$ is a Cartier 
 divisor. 

The relative ampleness of $P$ is an easy consequence of the fact that for 
each open affine subset $U$ of $S$, 
$\pi^{-1}(U)\setminus P = Spec({\cal A}(U))$ is affine.

Finally, one easily sees that $\pi_{*}({\cal O}_{C}(nP))$ can be identified 
with ${\cal A}^{(n)}$. So, the fact 
that $\bigcap_{n\geq 0} \pi_{*}{\cal O}_{C}(-nP) = (0)$ is a consequence of 
Lemma \ref{A0}. 
\QED
\freeline\\
Analogously, one can prove :
\begin{theorem}
$\widehat{\cal F}$ is a free $\widehat{\cal O}_{C}$-module, and 
the inclusion of $\cal W$ in \linebreak[4] $\Onegpower{S}{z}^{\oplus r}$ 
induces an 
isomorphism of sheaves of $
\widehat{\cal O}_{C}$ - modules:
$$\Phi : \widehat{\cal F} \stackrel{\sim}{\longrightarrow} 
\widehat{\cal O}_{C}^{\oplus r}.$$
 The intersection $\bigcap_{n\geq 0} \pi_{*}{\cal F}(-nP)$ vanishes. 
{}Furthermore, $F = \gamma(\pi_{*}{\cal F}) - 
\gamma(R^{1}\pi_{*}{\cal F})$.
\end{theorem}
The proofs of the first statements are  pure analogies to  Theorem 
\ref{Schur curve}. The very last assertion is a consequence of  
Theorem \ref{cohomology}.
\QED
\begin{definition}
{}For a given Schur pair $({\cal A},{\cal W})$ of rank $r$ and index $F$, 
we define 
$$\eta_{r,F}({\cal A},{\cal W}) :=(C,\pi,S,P,\rho,{\cal F},\Phi)$$
with the  objects described above. This defines a map
$$\eta_{r,F}:{\frak S}^{r}_{F}(S) \longrightarrow 
{\frak D}^{r}_{F}(S).$$
\end{definition}
Now we are ready to prove the converse of Theorem \ref{morph}.
\begin{theorem}
Homomorphisms of Schur pairs induce homomorphisms of the corresponding 
geometric data.
\end{theorem}
{\bf Proof}\freespace
Let $(\alpha,\xi):({\cal A}', {\cal W}')\rightarrow ({\cal A}, {\cal W})$ 
be a homomorphism of Schur pairs. We want to construct a homomorphism
$$(\alpha,\beta,\Psi):(C,\pi,S,P,\rho,{\cal F},\Phi)\rightarrow
(C',\pi',S',P',\rho',{\cal F}',\Phi')$$
of the associated geometric data. We proceed in two steps:

First we assume that $({\cal A}, {\cal W})=(\alpha^{(*)}{\cal A}', 
\alpha^{(*)}{\cal W}')$. This gives us a morphism

$$\beta: C=Proj(grd(\alpha^{(*)}{\cal A}'))\rightarrow C'$$
which makes the following diagram commute 
$$
\begin{array}{rcl}
C & \stackrel{\beta}{\rightarrow} & C'\\
\pi \downarrow &&\downarrow \pi'\\
S & \stackrel{\alpha}{\rightarrow} & S'.
\end{array}
$$
Note that, in general, $C$ is different from the fibre product $C'\times_{S'} S$ 
(cf. Section \ref{base change}). 

Since $C\setminus P=Spec(\alpha^{(*)}{\cal A}')$ maps to   
 $Spec({\cal A}')  = C'\setminus P'$, we get $\beta^{-1}(P') = P$.  Then, 
from the construction of $P$ in
  Theorem \ref{Schur curve}, it is also clear that $\beta^{*}P'=P$.

Recall that 
the local trivializations $\rho$ and $\rho'$ have been defined by the 
inclusions 
$\alpha^{(*)}{\cal A}' \subseteq \Onegpower{S}{z}$ and 
${\cal A}' \subseteq \Onegpower{S'}{z}$. So it is obvious that 
$\rho=\widehat{\beta}^{*}(\rho')$.

Finally, we consider the sheaves $\cal F$ and ${\cal F}'$. 
As ${\cal F}= (grd(\alpha^{(*)}{\cal W}'))^{\sim}$ and  
${\cal F}'= (grd({\cal W}'))^{\sim}$,  there is a natural map 
$\beta^{*}{\cal F}'\rightarrow {\cal F}$
 which is an isomorphism near $P$.
\freeline

Now let us return to the general case. As  $({\cal A}', {\cal W}')$ 
is a Schur pair, the induced object $(\alpha^{(*)}{\cal A}', \alpha^{(*)}{\cal W}')$ is 
also a Schur pair, i.e., the given homomorphism decomposes as follows:
$$({\cal A}', {\cal W}') \stackrel{(\alpha,id)}{\longrightarrow} 
(\alpha^{(*)}{\cal A}', \alpha^{(*)}{\cal W}')
\stackrel{(id_{S},\xi)}{\longrightarrow} ({\cal A}, {\cal W}).$$
So it just remains to consider the case where $S=S'$, $\alpha=id_{S}$, 
${\cal A}'\subseteq {\cal A}$ and $\xi({\cal W}')\subseteq {\cal W}$. 

The inclusion ${\cal A}'\subseteq {\cal A}$ induces 
$grd({\cal A}')\hookrightarrow grd({\cal A})$ and, therefore, a morphism
$$\beta: Proj(grd({\cal A})) = C \rightarrow C' = Proj(grd({\cal A}'))$$ 
which restricts to
$$\beta : Spec({\cal A}) \rightarrow Spec({\cal A}')$$
and which, in addition, is an isomorphism near $P$. 
Therefore, $\beta$ fits into the diagram
$$
\begin{array}{ccccc}
P \subset & C & \stackrel{\beta}{\rightarrow} & C'&\supset P'\\
\sim\searrow&\pi \downarrow &&\downarrow \pi'&\swarrow\sim\\
 &S&=& S
\end{array}
$$
and $\beta^{*}(P') = P$. The statement $\rho=\widehat{\beta}^{*}(\rho')$ is obvious. 

Now let us define the homomorphism of sheaves. We know that 
 $\xi({\cal W}') 
\subseteq {\cal W}$. Since $({\cal A},{\cal W})$ is a Schur pair, 
this implies:
$$\xi({\cal A}\cdot{\cal W}') ={\cal A}\cdot\xi({\cal W}')
\subseteq {\cal W}.$$
$\xi$ is determined by the images of the basis elements 
$e_{1},\ldots,e_{r'}$. By the definition of $\xi$,  
either $ord(\xi(e_{j}))\leq 0$ or $\xi(e_{j})= 0$.
 Therefore, $\xi(({\cal A}\cdot{\cal W}')^{(n)})\subseteq {\cal W}^{(n)}$,
 for all $n\in \Bbb{Z}$, and
 $\xi$ induces a homomorphism $\xi:grd({\cal A}\cdot {\cal W}') 
\longrightarrow grd({\cal W})$, i.e., a homomorphism of sheaves
$$\Psi : grd({\cal A}\cdot {\cal W}')^{\sim} = \beta^{*}{\cal F}'
\longrightarrow {\cal F} = grd({\cal W})^{\sim}.$$
Obviously, for this homomorphism $\Psi$, $\xi$ is recovered by $\xi = 
\pi_{*}((\rho\circ\Phi)\circ 
\widehat{\Psi}\circ \widehat{\beta}^{*}(\rho'\circ \Phi')^{-1})$.
\QED
\begin{corollary}
$\eta$ is a contravariant functor from the category of Schur pairs 
to the category of geometric data.
\QED
\end{corollary}
\begin{theorem}
The Krichever functor  $\chi$ and the functor $\eta$ are equivalences of the 
categories $\frak D$ and 
$\frak S$ and inverse to each other. Under this categorical equivalence, the 
subcategory ${\frak D}'$ 
corresponds to ${\frak S}'$.
\end{theorem}
{\bf Proof}\freespace
This is an easy consequence of Theorem II.5.14 \cite{H1}.
\QED
\freeline

\section{Applications}
\label{APPL}
Once the correspondence between geometric data and Schur pairs is 
established,  we are, of course, interested in seeing how this relation 
works practically. For example, assume that the given geometric 
objects have additional properties. How do these properties display in the 
corresponding Schur pair?

On the other hand, we had to impose some strong conditions on our 
family of curves and the sheaf on it (cf. Definitions \ref{total C} 
and \ref{sheaf F}). How substantial are these conditions? Are there still 
interesting and significant examples?

\subsection{Translation of geometric properties}
\label{geom. properties}

Let $(C,\pi,S,P,\rho,{\cal F},\Phi)$ be a geometric datum of rank $r$ and 
index $F$,  and $({\cal A},{\cal W})$ the associated Schur pair.
 A particular question is: What happens if $C$ or $\cal F$ is $S$-flat?
\begin{lemma}
\label{flat}
The sheaf $\cal F$ is flat over $S$ if and only if ${\cal W} \subset 
\Onegpower{S}{z}^{\oplus r}$ is locally free. $\pi$ is a flat morphism 
if and only if ${\cal A} \subset \Onegpower{S}{z}$ is locally free.
\end{lemma}
{\bf Proof}\freespace
Flatness is a local property. So we may assume that $S=Spec(R)$, $R$ a 
noetherian ring. 
We know that ${\cal O}_{C}(P)$ is ample on $C$ relative to $S$. Let 
$N$ be so that ${\cal O}_{C}(N\cdot P)$ is very ample relative to $S$. 
We know from the proof of \cite{H1}, Thm. III.9.9., that $\cal F$ is $S$ 
- flat if and only if, for sufficiently large $n$, $\pi_{*}({\cal F}
(nN\cdot P))$ is a locally free sheaf of ${\cal O}_{S}$ - modules of 
finite rank. Remember that $\Phi$ and $\rho$ induce an isomorphism of $\pi_{*}
({\cal F}(nN\cdot P))$ 
with ${\cal W}  \cap  \power{{\cal O}_{S}}{z}^{\oplus r}\cdot z^{-nN}$ 
(cf. Corollary \ref{order}). By assumption, ${\cal W} \cap 
\spower{{\cal O}_{S}}{z}{r}$ is coherent, hence is of finite rank. As $S$ 
is assumed to be noetherian, this implies that, for all $m$, ${\cal W}  
\cap  \power{{\cal O}_{S}}{z}^{\oplus r}\cdot z^{-m}$ is also of finite 
rank. This implies that the $S$ - flatness of $\cal F$ is equivalent to the 
local freeness of ${\cal W}  \cap  \power{{\cal O}_{S}}{z}^{\oplus r}
\cdot z^{-nN}$ for sufficiently large $n$. 
But on the other side, we also know that $\snegpower{{\cal O}_{S}}{z}{r}/
({\cal W}+\spower{{\cal O}_{S}}{z}{r})$ is coherent. 

Set $W:=H^{0}({\cal W})$. Then 
$\negpower{R}{z}^{\oplus r}/(W+\power{R}{z}^{\oplus r})$ is a finitely 
generated
$R$-module, hence, for sufficiently large $n$:
$$W+\power{R}{z}^{\oplus r}\cdot z^{-nN} = \negpower{R}{z}^{\oplus r}.$$

Denote by $\{e_{1},\ldots,e_{r}\}$ the standard basis of the $\power{R}{z}$ 
- module $\power{R}{z}^{\oplus r}$. Then, for all $i=1,\ldots,r$ and  $j > nN$, 
there are elements
$$w_{i,j} = e_{i}\cdot z^{-j} + \textrm{ terms of lower order } \in W.$$
Let $\bar{W}$ be the free $R$-submodule of $W$ generated by these elements. 
Then, of course,
$$W = \bar{W} \oplus (W\cap \power{R}{z}^{\oplus r}\cdot z^{-nN}),$$
and we see that $W$ is locally free if and only if this holds true for 
$W\cap \power{R}{z}^{\oplus r}\cdot z^{-nN}$. 

The proof of the second statement may be completed in the same way.
\QED
\freeline

At this time, let us outline a result which follows immediately from Lemma \ref{A0}. 
\begin{lemma}
\label{downtwist}
If $(C,\pi,S,P,\rho,{\cal F},\Phi)$ is a geometric datum then 
$\pi_{*}{\cal O}_{C}(-nP)$ and $\pi_{*}{\cal F}(-nP)$  vanish  for sufficiently 
large $n\in\Bbb{N}$.
\QED
\end{lemma}

Now we turn our attention to the stability of sheaves. 
\begin{definition}
We call  $\cal F$  {\em strongly semistable with respect to the section $P$}
iff 
there is an integer $N$ such that
$$\pi_{*}{\cal F}(N\cdot P) = R^{1}\pi_{*}{\cal F}(N\cdot P) = 0.$$
\end{definition}

Later on we will see that this notion of semistability is the most convenient 
one for the examination of commutative algebras of differential operators 
corresponding to sheaves over relative curves.

Translating the last definition in terms of Schur pairs, we get immediately:
\begin{lemma}\label{stab.Schur}
$\cal F$ is strongly semistable with respect to $P$ if and only if
$${\cal W} \oplus \Opower{S}{z}^{\oplus r}\cdot z^{-N} = \Onegpower{S}{z}^
{\oplus r},$$
for some $N\in \Bbb{Z}$.
\QED
\end{lemma}
\begin{corollary}
If $\cal F$ is strongly semistable with respect to $P$ then, in particular, 
$\cal F$ is flat over $S$.
\QED
\end{corollary}

\freeline

\begin{definition}
A coherent sheaf $\cal F$ on $C$ is called {\em simple} if 
$$ \pi_{*} {\cal E}nd_{{\cal O}_{C}}({\cal F}) = {\cal O}_{S}.$$
\end{definition}

In the set-up of Schur pairs it is easier to handle isomorphisms than 
 homomorphisms. That is why we are 
interested in the following statement. \begin{lemma}
\label{simple}
Let $S$ be reduced and 
assume that, for each point $s\in S$, the residue field $k(s)$ is infinite. Then 
$\cal F$ is simple if and
 only if  $ \pi_{*} {\cal A}ut_{{\cal O}_{C}}({\cal F}) = {\cal O}_{S}^{*}$.
\end{lemma}
{\bf Proof}\freespace
One implication is obvious. To prove the other one, we can assume, without
loss of generality, that
 $S = Spec(R)$ is affine. We want to prove 
$$End_{{\cal O}_{C}}({\cal F}) = R$$
under the assumption that $Aut_{{\cal O}_{C}}({\cal F}) = R^{*}$. Obviously,  
$R$ is contained in
 $End_{{\cal O}_{C}}({\cal F})$. Now assume that $End_{{\cal O}_{C}}({\cal F})$ 
contains an element
  $\phi$ which does not belong to $R$. 
 For $r \in R$, we consider the endomorphism $r+\phi$ and restrict it to the 
fibres of $\pi$: 
 $(r+\phi)_{s} := (r+\phi)|_{C_{s}}$. Since $S$ is reduced, $r+\phi$ is an
 isomorphism if and only if, 
 for all $s\in S$, $(r+\phi)_{s}$ is an isomorphism. We define
$$S(r) := \{s\in S/ (r+\phi)_{s} \textrm{ is an isomorphism }\}.$$
Obviously, these are open, possibly empty, subsets of $S$. Now we show that, 
for each $s\in S$, there
 is an element $r\in R$ such that  $s \in S(r)$:\\
Let us fix $s\in S$. We write 
$${\cal G}_{r} := ker((r+\phi)_{s}) .$$
Then ${\cal G}_{r}$ is a subsheaf of ${\cal F}_{C_{s}}$. We show that the 
sheaves ${\cal G}_{r}$ and ${\cal G}_{r'}$ 
 intersect only in the zero section whenever $r-r'$ is not contained in the
 prime ideal defining $s$ in $Spec(R)$.
  This can be seen locally on $C_{s}$. Let $V$ be an open affine subset of 
$C_{s}$, and $F := {\cal F}(V)$. 

Assume  there is an element $f\in F$ such that 
$(r+\phi)_{s}(f) = (r'+\phi)_{s}(f)=0$. Then $(r-r')\cdot f = 0$, which 
implies  that $f=0$.

As $k(s)$ is assumed to be infinite, there are infinitely many elements 
$r\in R$ so that their pairwise differences are
 not contained in the ideal of $s$.

Hence we obtain an infinite chain
$${\cal G}_{1}\subseteq\ldots\subseteq {\cal G}_{1} \oplus 
{\cal G}_{2}\oplus \ldots \subseteq {\cal F}.$$
Since $C_{s}$ is a noetherian scheme and ${\cal F}_{C_{s}}$ is coherent, this
 chain must become stationary, i.e., 
there are  infinitely many $r\in R$ such that 
 $(r+\phi)_{s}$ is an injective homorphism between two sheaves with the same 
Hilbert polynomial, hence it must be 
 an isomorphism, and this means  $s \in S(r)$ for all those $r$. 

This proves that the sets $S(r)$ form a covering of $S$. Now, $(r+\phi)|S(r)$ 
is an isomorphism, hence corresponds 
to an element of $H^{0}(S(r),{\cal O}_{S_{r}})^{*}$. Therefore, 
$\phi|S(r) \in H^{0}(S(r),{\cal O}_{S_{r}})$, i.e., 
$\phi \in H^{0}(S,{\cal O}_{S})=R$. This is a contradiction. 
\QED
\freeline

Now let us return to Schur pairs. The sheaf of groups
$Isom_{\Opower{S}{z}}(\Opower{S}{z}^{\oplus r})$ acts on $\frak{G}^{r}_{F}(S)$ , 
for every $F\in K(S)$. Using Corollary
 \ref{ident. Schur} we draw the following two conclusions:
\begin{corollary}
Let $S$ be as in Lemma \ref{simple}. Then ${\cal W}$ corresponds to a simple 
sheaf if and only if
$$Stab_{{\cal W}} Isom_{\Opower{S}{z}}(\Opower{S}{z}^{\oplus r}) = 
{\cal O}_{S}^{*}.$$
\QED
\end{corollary}
\begin{corollary}
Again let $S$ be as above. Then a sheaf ${\cal F}$ belonging to a geometric 
datum $(C,\pi,S,P,\rho,{\cal F},\Phi)$ is 
simple if and only if for all equivalent geometric data
 $(C,\pi,S,P,\rho,{\cal F},\Phi')$:
$$\Phi' = \lambda \Phi$$
for some $\lambda\in H^{0}(S,{\cal O}_{S})^{*}$.
\QED
\end{corollary}
\freeline

Now let us see to what determinant line bundles correspond. Assume that 
$S$ is noetherian, regular and 
separated.  ${\cal W}$ is a quasicoherent subsheaf of $\Onegpower{S}{z}^
{\oplus r}$. We know that 
$\Onegpower{S}{z}^{\oplus r}/({\cal W} + \Opower{S}{z}^{\oplus r})$ is 
coherent. Together with the fact that the base scheme $S$ is noetherian, 
this implies that there is an integer $N$ satisfying
\begin{equation}
\label{determinant}
\Onegpower{S}{z}^{\oplus r}={\cal W} + \Opower{S}{z}^{\oplus r}\cdot z^{-N}.
\end{equation}
Therefore,
$$
\begin{array}{rcl}
{\cal W}/{\cal W}^{(N)} & = & 
({\cal W} + \Opower{S}{z}^{\oplus r}\cdot z^{-N})/
\Opower{S}{z}^{\oplus r}\cdot z^{-N} \\
& = & \Onegpower{S}{z}^{\oplus r}/\Opower{S}{z}^{\oplus r}\cdot z^{-N}
\end{array}
$$
which is a trivial sheaf of ${\cal O}_{S}$ - modules. Subsequently,  
$det({\cal W})=det({\cal W}^{(N)})$ is a well-defined line bundle on $S$. 
Note that this definition does not depend on the choice of the integer $N$ 
occuring in the condition (\ref{determinant}). Furthermore, additivity holds 
for exact sequences.

On the other hand, for the given sheaf $\cal F$ of ${\cal O}_{C}$-modules  we 
may consider the so-called {\em determinant of the cohomology} (after P.~Deligne)
$$\lambda({\cal F}) := det(\pi_{*}{\cal F})\otimes (det (R^{1}
\pi_{*}{\cal F}))^{-1}.$$
We can prove the following result. 
\begin{proposition}
$\lambda({\cal F})\cong det({\cal W})$.
\end{proposition}
{\bf Proof}\freespace
This is an easy consequence of Proposition \ref{cohomology}. Using 
the fact proven there, we get:
$$\lambda({\cal F}) \cong det ({\cal W}^{(0)}) \otimes det 
(\Onegpower{S}{z}^{\oplus r}/({\cal W} + \Opower{S}{z}^{\oplus r}))^{-1}.$$
Now we consider the exact sequences of quasicoherent sheaves on $S$:
$$0\rightarrow {\cal W}/{\cal W}^{(0)}\rightarrow 
\Onegpower{S}{z}^{\oplus r}/\Opower{S}{z}^{\oplus r}
\rightarrow \Onegpower{S}{z}^{\oplus r}/({\cal W} + \Opower{S}{z}^{\oplus r})
\rightarrow 0$$
and 
$$0\rightarrow {\cal W}^{(0)} \rightarrow {\cal W} 
\rightarrow {\cal W}/{\cal W}^{(0)} \rightarrow 0.$$
The statement of the lemma follows then from the additivity of $det$.
\QED
\freeline\\
{\bf Remark}\freespace
In the case that $\cal F$ is flat and $S$ is separable, $det({\cal W})$ is 
again well-defined, even if $S$ is not regular. Consequently this determinant 
generalizes the determinant of the cohomology. 
\freeline

\subsection{Examples of geometric data}
\label{Eogd}

First we prove a criterion which will be highly useful for the 
construction of examples.
\begin{proposition}
\label{lifting1}
Assume that the base scheme $S$ satisfies: $H^{1}(S,{\cal O}_{S}) = 0$, 
and that, for the section $P$, ${\cal I}_{P}/{\cal I}_{P}^{2}$ is free 
of rank 1.
Then 
\begin{enumerate}
\item $\widehat{\cal O}_{C} \cong \Opower{P}{z}$.
\item If $\cal F$ is a  coherent sheaf of ${\cal O}_{C}$ 
- modules
such that
\begin{itemize}
\item $\cal F$ is locally free in a neighborhood of  $P$,
\item ${\cal F}|P \cong {\cal O}_{P}^{\oplus r}$
\end{itemize}
then $\widehat{\cal F} \cong \widehat{\cal O}_{C}^{\oplus r}$.
\end{enumerate}
\end{proposition}
{\bf Remark}\freespace
For example, the cohomological condition is fulfilled for all affine schemes $S$.
\freeline\\
{\bf Proof of the proposition}\freespace
At first, observe that the condition that \linebreak[4] $H^{1}(S,{\cal O}_{S})$
vanishes is, of course, equivalent to $H^{1}(P,{\cal O}_{P}) = 0$. 

By Lemma \ref{powerseries} and the remark following  it the first claim is equivalent to:
$$H^{0}(C,{\cal I}/{\cal I}^{n}) \surjection 
H^{0}(C,{\cal I}/{\cal I}^{2}) 
\quad \forall n \in \Bbb{N}, n\geq 2.$$
One easily sees that this is the case if and only if for all 
$n\geq 2$:
$$H^{0}(C,{\cal I}/{\cal I}^{n+1}) \surjection 
H^{0}(C,{\cal I}/{\cal I}^{n}).$$
We have the exact sequence of sheaves of ${\cal O}_{P}$-modules
\begin{equation}
\label{equ.1}
0 \rightarrow {\cal I}^{n}/{\cal I}^{n+1} \rightarrow 
{\cal I}/{\cal I}^{n} \rightarrow {\cal I}/{\cal I}^{n+1}\rightarrow 0
\end{equation}
which induces a long exact sequence of cohomology groups
$$0 \rightarrow H^{0}({\cal I}^{n}/{\cal I}^{n+1}) \rightarrow 
H^{0}({\cal I}/{\cal I}^{n}) \rightarrow H^{0}({\cal I}/{\cal I}^{n+1}) 
\rightarrow H^{1}({\cal I}^{n}/{\cal I}^{n+1}) \rightarrow \ldots$$
By assumption, ${\cal I}^{n}/{\cal I}^{n+1} =  ({\cal I}/{\cal I}^{2})^{n} \cong
 {\cal O}_{P}$. So our assumption on  $S$ implies 
that  $H^{1}({\cal I}^{n}/{\cal I}^{n+1}) = 0$ for all $n\in \Bbb{N}$ and we are done.

Now we come to the second part. By Lemma \ref{COMP.F} the claim is equivalent to the 
following fact:
\begin{equation}
\label{equ.0}
H^{0}({\cal F}/({\cal I}^{n}\otimes {\cal F})) \surjection 
H^{0}({\cal F}/({\cal I}\otimes {\cal F})) \quad, \forall n\in \Bbb{N}.
\end{equation}
We consider one more exact sequence of coherent sheaves of 
${\cal O}_{C}$ - modules:
\begin{equation}
\label{equ.2}
0 \rightarrow {\cal I}^{n} \rightarrow 
{\cal O}_{C} \rightarrow {\cal O}_{C}/{\cal I}^{n} \rightarrow 0.
\end{equation}
Since $\cal F$ is locally free in some neighborhood of $P$, and 
${\cal O}_{C}/{\cal I}^{n} = 0$ outside $P$, the sequences 
(\ref{equ.1}) and (\ref{equ.2}) stay exact when we tensor with $\cal F$. 
So we get
\begin{equation}
\label{equ.3}
0 \rightarrow ({\cal I}^{n-1}/{\cal I}^{n})\otimes {\cal F} \rightarrow 
({\cal O}_{C}/{\cal I}^{n})\otimes {\cal F} \rightarrow 
({\cal O}_{C}/{\cal I}^{n-1})\otimes {\cal F} \rightarrow 0
\end{equation}
and (\ref{equ.2}) implies: 
$({\cal O}_{C}/{\cal I}^{n})\otimes {\cal F} \cong 
{\cal F}/({\cal I}^{n}\otimes {\cal F})$. Now we write down the 
long exact sequence of cohomology groups induced by (\ref{equ.3}):
\begin{equation}
\label{equ.4}
\begin{array}{cccccc}
0 & \rightarrow & H^{0}(({\cal I}^{n-1}/{\cal I}^{n})\otimes {\cal F}) & 
\rightarrow &
 H^{0}({\cal F}/({\cal I}^{n}\otimes {\cal F})) & \rightarrow \\
&\rightarrow & H^{0}({\cal F}/({\cal I}^{n-1}\otimes {\cal F})) & \rightarrow & 
H^{1}(({\cal I}^{n-1}/{\cal I}^{n})\otimes {\cal F}).
\end{array}
\end{equation}
Since the restriction of $\cal F$ to $P$ is free, and ${\cal I}^{n-1}/{\cal I}^{n}$ 
is isomorphic to 
${\cal O}_{P}$, $({\cal I}^{n-1}/{\cal I}^{n})\otimes_{{\cal O}_{C}} {\cal F}$ is a 
free ${\cal O}_{P}$-module. 
 Therefore we finally get:
$$H^{1}(({\cal I}^{n-1}/{\cal I}^{n})\otimes {\cal F})
\cong H^{1}(P,{\cal O}_{P})^{\oplus r}  = 0,$$
which, together with the sequence (\ref{equ.4}), implies the 
surjectivity in (\ref{equ.0}).
\QED
\freeline

Now we come to explicit examples.
\subsubsection{Trivial families of curves}
The easiest case, but which is not without interest, is the one of a 
trivial family of curves with some  sheaf on it.
\begin{proposition}
\label{CxS}
Let $K$ be a complete, integral curve over some field $k$ and $p\in K$ 
a smooth, $k$-rational point. Let $S$ be a locally noetherian 
$k$-scheme and set $C:= K\times_{Spec(k)} S$. Denote by $\pi$ the projection from
 $C$ to $S$, and by $P$ the section 
$\{ (p,s)/s\in S\} $. Then $\widehat{\cal O}_{C} \cong \Opower{P}{z}$.
\end{proposition}
{\bf Proof}\freespace
Let ${\cal J} = {\cal J}_{p}$ be the sheaf of ideals defining $p$ 
in $K$. Since $p$ is a smooth point, $\cal J$ is generated by one 
element $z$ near $p$, and we get:
\begin{itemize}
\item ${\cal J}/{\cal J}^{2} = [z]\cdot ({\cal O}_{K}/{\cal J})$;
\item $\widehat{\cal O}_{K} \cong \power{k}{z}$.
\end{itemize}
Since $C=K\times_{Spec(k)} S$ and $P=\{ p \} \times_{Spec(k)} S$, 
we obtain:
\begin{itemize}
\item ${\cal I}_{P} = {\cal J}\otimes_{k} {\cal O}_{S}$, hence 
${\cal I}_{P}/{\cal I}_{P}^{2} = 
[z]\cdot({\cal O}_{C}/{\cal I}_{P})=[z]\cdot{\cal O}_{P}$.
\item $\widehat{\cal O}_{C} = \widehat{\cal O}_{K}
\otimes_{k} {\cal O}_{S}$, i.e., 
$\widehat{\cal O}_{C} \cong \Opower{P}{z}$.
\end{itemize}
\QED
\freeline

\subsubsection{Elliptic curves}
\label{Elliptic curves}

To describe a nontrivial family of integral curves which fits into our set-up,  
we define a family of elliptic curves over $\Bbb{A}^{2}=\Bbb{A}^{2}_{\Bbb{C}}$ as follows:
\begin{equation}
\begin{array}{ccc}
C\textrm{ }&  := & \{(A,B,z_{0}:z_{1}:z_{2}) \in \Bbb{A}^{2}\times 
\Bbb{P}^{2} /
z_{0}z_{2}^{2}= z_{1}^{3} + A z_{0}^{2}z_{1} + B z_{0}^{3} \}\\
\downarrow  \pi & {} &{}\\
\Bbb{A}^{2}\textrm{ }&{}&{}
\end{array}
\end{equation}
A section of $\pi$ can be defined by  
$P:=\{(A,B,0:0:1)/(A,B)\in\Bbb{A}^{2}\}$. One 
easily sees that $\pi$ is a flat, projective morphism with reduced, 
irreducible fibres of dimension 1 and that $C$ is reduced.

We want to study the conormal sheaf of $P$. Since 
$P$ does not intersect the hyperplane 
$\Bbb{A}^{2}\times(z_{2}=0)$, we can restrict our consideration to 
its affine complement. Let us denote $y_{i} := z_{i}/z_{2}$ for 
$i=0,1$. Then $C\cap(z_{2}\neq 0)\subset \Bbb{A}^{2}\times 
\Bbb{A}^{2}$ is given by the equation
\begin{equation}
\label{elliptic}
y_{0} = y_{1}^{3} + A y_{0}^{2}y_{1} + B y_{0}^{3}.
\end{equation}
Let $R$ be the affine coordinate ring of $C\cap(z_{2}\neq 0)$.
The ideal $I$ of $P$ in $R$ is generated by $y_{0}$ 
and $y_{1}$. But $y_{0} \in I^{2}$, i.e. $I/I^{2} = y_{1}R/I$. 
This implies that
 ${\cal I}/{\cal I}^{2}$ is a trivial line bundle. 
\freeline\\

Note that $\widehat{\cal O}_{C} \cong \Opower{P}{z}$ since $H^{1}(\Bbb{A}^{2}) =0$.

Now we want to find a suitable (formal) local parameter on $C$ along 
$P$. From the above calculation we know that 
$y_{1}=z_{1}/z_{2}$ is such a local parameter, i.e., the formal 
completion $\widehat{R}$ of $R$ with respect to $I$ equals  
$k\/[\/A,B\/]\/[[\/y_{1}\/]]$. We claim that for $\alpha := 
\sqrt{z_{0}/z_{1}}$, $\widehat{R} \cong 
k\/[\/A,B\/]\/[[\/\alpha\/]]$. We use the equation (\ref{elliptic}) and 
calculate:
$$
\begin{array}{rcl}
y_{0} & = & y_{1}^{3} + A y_{0}^{2}y_{1} + B y_{0}^{3}\\
y_{0}/y_{1} & = & y_{1}^{2} + A y_{0}^{2} + B y_{0}^{3}/y_{1}\\
-y_{1}^{2} - A y_{0}^{2} & = & (B y_{0}^{2} -1)y_{0}/y_{1}.
\end{array}
$$
$(B y_{0}^{2} -1)$ is an invertible element of 
$k\/[\/A,B\/]\/[[\/y_{1}\/]]$. Therefore 
$$z_{0}/z_{1}= y_{0}/y_{1} = 
(B y_{0}^{2} -1)^{-1}(-y_{1}^{2} - A y_{0}^{2}) \in 
k\/[\/A,B\/]\/[[\/y_{1}\/]]$$
is an element of order $-2$ with leading coefficient 1. So, 
the square root of $z_{0}/z_{1}$ is a well-defined monic element 
$\alpha$ of $\widehat{R}$ of order $-1$. This implies that 
$\widehat{R}  \cong k\/[\/A,B\/]\/[[\/\alpha\/]]$.

Now we construct the corresponding subring $\cal A$ of 
$k\/[\/A,B\/]\/[[\/\alpha\/]]$. 
We see that $P$ is exactly the intersection of 
$C$ with the hyperplane $(z_{0}=0)$. Therefore, the affine ring 
of coordinates of $C\setminus P$ is 
$$k\/[\/A,B,z_{1}/z_{0}, z_{2}/z_{0}\/]/((z_{2}/z_{0})^{2}- 
(z_{1}/z_{0})^{3} -A (z_{1}/z_{0})^{2} - B).$$

We express the generating elements of this $k\/[\/A,B\/]$ - algebra 
in terms of the above chosen formal parameter $\alpha$: 
$$
\begin{array}{rcl}
z_{1}/z_{0} & = & \alpha^{-2}\\ 
z_{2}/z_{0} & = & \alpha^{-2}\cdot y_{1}^{-1}.
\end{array}
$$
$y_{1}^{-1}$ is an element of $k\/[\/A,B\/]\/[[\/\alpha\/]]$. We want
 to find out its special form. Using the equation (\ref{elliptic}) we get
$$
\begin{array}{rcl}
1/y_{1}^{2} & = & y_{1}/y_{0} + A (y_{0}/y_{1})+ B (y_{0}/y_{1})^{2}\\
{} & = & \alpha^{-2} + A \alpha^{2}+ B\alpha^{4}.
\end{array}$$
So, finally:
 $${\cal A } = k\/[\/A,B,\alpha^{-2}, \alpha^{-2}\cdot\sqrt{\alpha^{-2}
 + A \alpha^{2}+ B\alpha^{4}}\/].$$
The term $\alpha^{-2}\in {\cal A}$ reflects the 2:1 covering 
$$C\rightarrow \Bbb{A}^{2}\times \Bbb{P}^{1}.$$
\subsubsection{Families of line bundles over a curve}
\label{Poincare}
Let $C$ be a  complete integral complex curve  and 
$p\in C$ a  point. We choose a formal local 
trivialization $\rho$ of $C$ near $p$ and construct the corresponding 
subring $A\subset \negpower{\Bbb{C}}{z}$. 

We take the Picard variety $Pic^{n}(C)$ of $C$, for some $n$,  as a base scheme. As 
described in Proposition \ref{CxS},
 $\rho$ extends to a local 
trivialization of $C\times Pic^{n}(C)$ near $\{p\}\times Pic^{n}(C)$ 
and we get for the corresponding sheaf of ${\cal O}_{Pic^{n}(C)}$-
algebras:
$${\cal A} = A\otimes_{\Bbb{C}}{\cal O}_{Pic^{n}(C)}\subseteq 
\Onegpower{Pic^{n}(C)}{z}.$$

Now we consider the Poincar\'{e} bundle ${\cal P}^{n}_{C}$ of degree 
$n$ on $C$ (normalized with respect to the fixed point $p$). 
${\cal P}^{n}_{C}$ is a line bundle on  $C\times  Pic^{n}(C)$ 
satisfying:
\begin{itemize}
\item ${\cal P}^{n}_{C}|C\times \{L\} \cong L$ for every 
$L\in Pic^{n}(C)$,
\item ${\cal P}^{n}_{C}|\{p\}\times Pic^{n}(C)$ is trivial,
\item ${\cal P}^{n}_{C}$ is flat over $Pic^{n}(C)$.
\end{itemize}
For more details see \cite{LB}.

Let $U\subset Pic^{n}(C)$ be an open affine subset. We apply  Proposition \ref{lifting1} 
and conclude that 
${\cal P}^{n}_{C}|(C\times U)$ satisfies Condition 1 of Definition \ref{sheaf F}. So we can 
construct (for some local trivialization) the corresponding sheaf 
${\cal W}(U)\subseteq \Onegpower{U}{z}$ and we obtain a Schur pair $({\cal A}(U), {\cal W}(U))$. 
(Keep in mind that
 $\cal A$ is globally defined whereas $\cal W$ is not!)

${\cal W}(U)$ defines a map from $U$ to 
${\frak G}^{1}_{\mu}(Spec(\Bbb{C}))$, $\mu = n+1-g(C)$, and as the 
image we get $U$ as a subset of the Grassmannian.

M.~Mulase proved in \cite{M1} that every finite dimensional integral 
manifold of the KP-flows on some quotient of the Grassmannian 
${\frak G}^{1}_{\mu}(Spec(\Bbb{C}))$ has the linear structure of a Jacobian 
of a curve. The quotient has been taken in order to eliminate different 
local trivializations. In particular, the differences arising from the local construction 
of $\cal W$ cancel out completely.
  So, the above result implies that the integral manifold also  carries 
 the algebraic-geometric structure of the Jacobian.
\subsubsection{Base change}
\label{base change}

Here,  we want to investigate the behaviour of geometric data under base changes. In general, 
the pull-back of a
 geometric datum over $S$ under a base change $\alpha:S'\rightarrow S$ does not
 give rise to another geometric datum.
  However, let us start with a positive example. 
\begin{lemma}
Let $(C,\pi,S,P,\rho,{\cal F},\Phi)$ be a geometric datum of rank $r$ over $S$
 and $\alpha:S'\rightarrow S$ a flat morphism.
 The fibre product construction then defines a collection 
$(C',\pi',S',P',\rho',{\cal F}',\Phi')$. We claim that this collection
  forms a geometric datum of rank $r$ over $S'$ and that 
$(\alpha,\alpha',id)$ is a homomorphism of geometric data
$$(\alpha,\alpha',id): (C',\pi',S',P',\rho',{\cal F}',\Phi')
\rightarrow 
(C,\pi,S,P,\rho,{\cal F},\Phi),$$ 
where $\alpha'$ is defined as the fibre product morphism,
$$
\begin{array}{rcl}
{}C' & \stackrel{\alpha'}{\rightarrow}  & C\\
\pi'\downarrow && \downarrow \pi\\
{}S' & \stackrel{\alpha}{\rightarrow}& S.
\end{array}
$$
\end{lemma}
{\bf Proof}\freespace
Let us check the properties listed in Definition \ref{def data}. Some
 of them are easy to see. Of course, $C'$ is a scheme, 
$\pi':C'\rightarrow S'$ is a locally projective morphism,  $P'$ is a relatively 
ample Cartier divisor and 
${\cal I}'/{\cal I}'^{2} = \alpha'^{*}({\cal I}/{\cal I}^{2})$ is free of rank 1 
on $P'$. 

In Lemma \ref{powerseries} we saw that the condition $\widehat{\cal O}_{C} 
\cong \Opower{P}{z}$ is equivalent to the
 fact that the section 
$1\in H^{0}(C, {\cal I}/{\cal I}^{2})= H^{0}(P,{\cal O}_{P})$ lifts to a 
section of ${\cal I}/{\cal I}^{n}$ for all 
$n\geq 2$. Working through the diagram
$$
\begin{array}{ccccccc}
H^{0}({\cal I}/{\cal I}^{n}) & \surjection & H^{0}({\cal I}/{\cal I}^{2}) & 
= & H^{0}({\cal O}_{P}) & \ni & 1\\
\alpha^{*} \downarrow&& \alpha^{*}\downarrow\\
H^{0}({\cal I}'/{\cal I}'^{n}) & \surjection & H^{0}({\cal I}'/{\cal I}'^{2}) & 
= & H^{0}({\cal O}_{P'}) & \ni & 1
\end{array}
$$
we obtain: $\widehat{\cal O}_{C'} \cong \Opower{P'}{z}$.

Of course, $\rho$ and $\alpha$ induce an isomorphism
$$\rho': \widehat{\cal O}_{C'} \stackrel{\sim}{\rightarrow} \Opower{P'}{z}.$$

Now we turn our attention to the sheaf ${\cal F}'$. It is easy to see that
$${\cal F}'/{\cal I}'\otimes {\cal F}' \cong \alpha'^{*}({\cal F}/{\cal I}\otimes 
{\cal F}) \cong {\cal O}_{P'}^{\oplus r}$$
and that ${\cal F}'$ is locally free near $P'$.

In order  to prove that the completion of ${\cal F}'$ along $P'$ is free, we have 
to show that the generating sections
$e_{1},\ldots, e_{r}$ of  $H^{0}({\cal O}_{P'}^{\oplus r}) = 
H^{0}({\cal F}'/{\cal I}'\otimes {\cal F}')$ lift to sections
 of ${\cal F}'/{\cal I}'^{n}\otimes {\cal F}'$ for all $n\in\Bbb{N}$. 
This is done as above. 

Again, $\Phi$ and $\alpha$ induce an isomorphism of 
$\widehat{\cal O}_{C'}$-modules
$$\Phi' : \widehat{{\cal F}'} \stackrel{\sim}{\rightarrow} 
\widehat{\cal O}_{C'}^{\oplus r}.$$

Finally, note that for sufficiently large $n\in\Bbb{N}$, 
$$\pi_{*}{\cal O}_{C}(-nP) =0 \textrm{ and } \pi_{*}{\cal F}(-nP) =0$$
(cf. Lemma \ref{downtwist}). Consequently, using \cite{H1}, Thm.II.9.3.,
$$\pi'_{*}{\cal O}_{C'}(-nP') = \alpha^{*}\pi_{*}{\cal O}_{C}(-nP)=0 
\textrm{ and } \pi'_{*}{\cal F}'(-nP') =\alpha^{*}\pi_{*}{\cal F}(-nP)=0$$
for sufficiently large $n$. This completes the list of the properties which we 
had to check and
 $(C',\pi',S',P',\rho',{\cal F}',\Phi')$ is really a geometric datum. 

The fact that $(\alpha,\alpha', id)$ is a morphism of geometric data is 
straightforward. \QED
\freeline

As we mentioned at the beginning, the pull-back of a geometric datum over
 $S$ via a morphism $\alpha:S'\rightarrow S$ 
does not always define a new geometric datum. A typical situation for this
 is the restriction to one point of the base scheme $S$. 

Let us give an example in terms of Schur pairs.  As a base scheme we choose 
$S=\Bbb{A}^{1}_{k} = Spec (\poly{k}{t})$, and our Schur pair is given by $(A,A)$ for
$$A = \poly{k}{t} \oplus t\poly{k}{t}z^{-1} \oplus \poly{\poly{k}{t}}{z^{-1}}\cdot z^{-2}.$$
Now we consider the fibre of the corresponding projective curve $C$ over 
 the point $0\in S$. Its affine part outside 
the section $P$ is given by the ring
$$A_{0} = k \oplus (t\poly{k}{t}/t^{2}\poly{k}{t})z^{-1} \oplus 
\poly{k}{z^{-1}}\cdot z^{-2}.$$
This is a cuspidal curve with an embedded point. In particular, the fibre
 is not reduced and therefore cannot be a
 part of a geometric datum over $k$. 

However, $(A,A)$ in fact induces a Schur pair over $S'=Spec(k)$, namely (A',A') with
$$A' = Im(A\subset \negpower{\poly{k}{t}}{z} \surjection 
\negpower{(\poly{k}{t}/t\cdot\poly{k}{t})}{z})
= k\oplus \poly{k}{z^{-1}} z^{-2}.$$
This ring corresponds to the cuspidal curve, i.e., to the integral
 component of the fibre passing through the section $P$.
\newpage 
\section{Families of commutative algebras of differential operators}
\label{Fam. DO}
In \cite{M1} M.~Mulase used the  equivalence of the category of Schur pairs
 and the category of quintets for a complete classification of commutative 
algebras of 
ordinary differential operators with coefficients in $\power{k}{x}$. This 
leads us to the natural question of 
whether it is possible to extend these results to the relative case, at 
least in some special situations. It is hard to do this in the set-up of 
sheaves. Therefore we will restrict our observations to the case of an 
affine base scheme 
$S=Spec(R)$, where $R$ is a commutative noetherian $k$-algebra for some field
 $k$ of characteristic zero.

Before beginning let us fix a convention: Whenever in this chapter we speak about
 Schur pairs we mean elements
 of ${\frak S}'$ and all geometric data occuring here belong to ${\frak D}'$.

 Let us start  analyzing our objects in this special case.
\subsection{Schur pairs over affine base schemes}
\label{affine base}
\begin{definition}
Let $A$ be an $R$-subalgebra of $\negpower{R}{y}$, and 
$r\in\Bbb{N}$. $A$ is said to be an {\em algebra of pure rank $r$}, 
if 
\begin{enumerate}
\item $r=gcd(ord(a)/a\in A)$ and
\item There are monic elements $a$ and $b$ of positive order in $A$ 
such that $gcd(ord(a),ord(b))=r$.
\end{enumerate}
\end{definition}
Let us see what these properties imply:
\begin{lemma}
\label{pure}
Let $A\subseteq \negpower{R}{y}$ be an $R$-subalgebra of pure rank $r$. 
Then there is a monic element $z\in \power{R}{y}$ of order $-r$ such that 
\begin{itemize}
\item $A\subseteq \negpower{R}{z}$,
\item $\negpower{R}{z}/(A+\power{R}{z})$ is a finitely generated 
$R$-module.
\end{itemize}
\end{lemma}
{\bf Proof}\freespace
Choose monic elements $a$ and $b$ of $A$ of positive order such that 
$gcd(ord(a),ord(b))=r$. Then there are natural numbers $i$ and $j$
 such that 
$$r=i(ord(a))-j(ord(b)).$$
Define $z:= a^{-i}b^{j}$. Since the inverse of a monic element of 
$\negpower{R}{y}$ is again a well-defined element of 
$\negpower{R}{y}$, we have constructed  a monic element \linebreak[4] $z\in 
\power{R}{y}$ of order $-r$. Now let us prove that the localization 
of $A$ by $a$ is contained in $\negpower{R}{z}$:

Of course, $z\in A_{a}$.  We choose an element $v$ of $A_{a}$ and denote its 
order by $\alpha$. $v$ has the
 form $\frac{w}{a^{m}}$ for some elements $w\in A$ and  $m\in\Bbb{N}$. $r$ 
divides the orders of $w$ and $a$,
  therefore $r$ divides $\alpha$.  Since $z$ is 
monic of order $-r$, there is some $n\in\Bbb{Z}$ and $v_{0}\in R$, 
such that $v-v_{0}\cdot z^{n} \in A_{a}$ is an element of order 
less than $\alpha$. Now the assertion is proved inductively. 

In particular, this shows that $A$ itself is contained in 
$\negpower{R}{z}$.

For the second part see the proof of Lemma \ref{a and b}.
\QED
\freeline\\
{\bf Remark}\freespace
The converse 
of this lemma is true, as well:\\
Let $A\subset \negpower{R}{y}$ be an $R$-subalgebra satisfying
\begin{itemize}
\item $A\subset \negpower{R}{z}$ for some monic element $z$ of order 
$-r$;
\item $\negpower{R}{z}/(A+\power{R}{z})$ is a finitely generated 
$R$-module
\end{itemize}
then $A$ is an $R$-algebra of pure rank $r$.
\freeline\\
Now let us define
\begin{definition}
As an {\em embedded Schur pair of rank $r$, index $F$ and level 
$\alpha$ over $Spec(R)$} we denote a pair $(A,W)$ consisting of
\begin{itemize}
\item $A\subseteq \negpower{R}{y}$ an $R$-subalgebra of pure rank $r$ satisfying
 $A\cap\power{R}{y}=R$ ;
\item $W\subseteq \negpower{R}{y}$ with $W\in 
\frak{G}^{1}_{F,\alpha}(Spec(R))$
\end{itemize}
such that $A\cdot W \subseteq W$.

We write $\frak{E}_{\alpha}\frak{S}'^{r}_{F}(R)$ for the set of 
embedded Schur pairs of rank $r$, index $F$ and level $\alpha$ 
over $Spec(R)$.
\end{definition}
This is a natural generalization of the notion of Schur pairs 
introduced by M.~Mulase \cite{M1}. Now we want to see how these 
objects are related to the Schur pairs on $Spec(R)$ we have defined earlier.
\begin{proposition}
For all $\alpha\in\Bbb{Z}$ there is a canonical one-to-one 
correspondence between Schur pairs of rank $r$ and index $F$ and 
embedded Schur pairs $(A,W)$ of rank $r$, index $F$ and level 
$\alpha$ with the extra-condition that
$$A\subset \negpower{R}{y^{r}}.$$
\end{proposition}
{\bf Proof}\freespace
The method of the proof has been outlined already in the remark on 
page \pageref{Grass Mulase}.

Let us start with an embedded Schur pair $(A,W)$ of rank $r$, index $F$ 
and level $\alpha$ such that $A\subset \negpower{R}{y^{r}}$. Set $z:=y^{r}$. 
Then, by Lemma \ref{pure},  $A\subset\negpower{R}{z}$ is an 
element of ${\frak G}^{1}_{G}(Spec(R))$ for some $G\in K(Spec(R))$. Now we identify:
$$
\begin{array}{rcl}
\power{R}{y}\cdot y^{-\alpha} & = & \bigoplus_{i=-\alpha}^{-\alpha+r-1} 
\power{R}{y^{r}}\cdot y^{i}\\
&=& \power{R}{z}^{\oplus r}.
\end{array}
$$
This identification extends to an isomorphism of $\negpower{R}{y}$ with 
$\negpower{R}{z}^{\oplus r}$ and so we end up with  $W\subset 
\negpower{R}{z}^{\oplus r}$. Since $\power{R}{y}\cdot y^{-\alpha}$ 
 translates into $\power{R}{z}^{\oplus r}$, $W$ gives an element of 
${\frak G}^{r}_{F}(Spec(R))$.

That also clarifies the inverse construction. We formally set 
$z:=y^{r}$ and translate the data back using Lemma \ref{pure} 
and its converse.
\QED\freeline

\subsection{Formal pseudo-differential operators}
\label{pseudo DO}
We  saw that embedded Schur pairs are closely related to Schur pairs, while 
Schur pairs themselves correspond to geometric data via the Krichever functor. 
Now a natural question is how to identify embedded Schur pairs which lead 
to ``similar'' geometric data, where similar means that they differ only by a 
very special change of the local trivializations. This is done with the 
help of formal pseudo-differential operators. Furthermore, formal 
pseudo-differential operators will 
be the main tool for the classification of commutative algebras of 
differential operators.\freeline

Consider the ring $\power{R}{x}$ of formal power series in one variable 
with coefficients in $R$, and write $\partial :=\frac{d}{dx}$. 
$\partial$ acts on $\power{R}{x}$ by derivation:
$$\partial(\sum_{i\geq 0}a_{i} x^{i}) = \sum_{i\geq 1} ia_{i}x^{i-1},$$
while, for any $n\in\Bbb{N}$, $\partial^{n}$ acts by repeated derivation:
$$\partial^{n} (f) := \partial(\partial^{n-1} f)$$
for $f\in\power{R}{x}$. $\partial^{0}$ is defined to be the identity.

For given elements $f,g\in\power{R}{x}$ and $n\in\Bbb{N}$ we define:
$$\begin{array}{rcl}
(f\partial^{n})(g) &=&f\partial^{n}(g),\\
(\partial^{n}f)(g)&=&\partial^{n}(fg).
\end{array}$$
In this way, the ring of ordinary differential operators with coefficients 
in $\power{R}{x}$, 
$D := \poly{\power{R}{x}}{\partial}$, turns out to be a subring of the 
endomorphism ring $End_{R}(\power{R}{x})$. 
The multiplication of elements of $D$ is determined by the {\em Leibniz rule}:
\freeline\\
For $f,g \in \power{R}{x}$ and $n\in \Bbb{N}$:
\begin{equation}
\label{Leib}
\partial^{n}(fg) = \sum_{i=0}^{\infty} {n \choose i} f^{(i)}
\partial^{n-i}(g).
\end{equation}
It is our aim to make the operator $\partial$ invertible. In fact, we want 
to introduce $\partial^{-1}$ with 
$\partial\partial^{-1} = \partial^{-1}\partial = 1$ and define a
 multiplication on the set
$$E:=\{\sum_{n\in \Bbb{Z}} f_{n}\partial^{n}/ f_{n}\in \power{R}{x}, 
f_{n} = 0 \textrm{ for } n \gg 0\}\supset D,$$
which is compatible with the multiplication on $D$. One can define $\partial^{-1}$ as an endomorphism on
 $\power{R}{x}$ by formal integration:
$$\partial^{-1}(\sum_{i\geq 0}a_{i} x^{i}) = \sum_{i\geq 0} 
\frac{a_{i}}{i+1}x^{i+1}.$$
But obviously, the so-defined $\partial^{-1}$ is not inverse to $\partial$ 
as an endomorphism of $\power{R}{x}$. 
That is why we consider the action of $\partial$ and $\partial^{-1}$ on the 
quotient $R$-module $\power{R}{x}/\poly{R}{x}$.
 One easily sees that one may interpret $\poly{\poly{R}{x}}{\partial, 
\partial^{-1}}$ as a subring of
  $End_{R}(\power{R}{x}/\poly{R}{x})$ and that for this embedding $\partial^{-1}$ 
is the inverse of $\partial$.

Observe that the Leibniz rule (\ref{Leib}) also holds true for negative $n$,
 i.e., for  formal integration, where,
 for arbitrary $n\in\Bbb{Z}$ and $i\in\Bbb{N}$, the binomial coefficient 
${n \choose i}$ is defined as follows:
$${n \choose i} := 
\frac{n\cdot(n-1)\cdot\ldots\cdot(n-i+1)}{i\cdot (i-1)\cdot\ldots\cdot 1}
\in \Bbb{N}.$$
For negative $n\in\Bbb{Z}$, the summation in (\ref{Leib}) is really an infinite 
one, while, for nonnegative $n$, it is finite.

The formula (\ref{Leib}) defines a  multiplication rule for elements of $E$ of the form
$$\sum_{n=-M}^{N} (\sum_{i=0}^{\alpha_{n}} f_{i,n} x^{i})\partial^{n}.$$
But now it is clear that the so-defined multiplication extends to a 
multiplication on all of $E$, which restricts
 to the usual one on $D$ and has the form:
$$
\begin{array}{l}
(\sum_{m=0}^{\infty}a_{m}\partial^{M-m})\cdot
(\sum_{n=0}^{\infty}b_{n}\partial^{N-n})=\\
  \hspace{3cm}=\sum_{m=0}^{\infty}\sum_{n=0}^{\infty}\sum_{i=0}^{\infty} {M-m \choose i}
a_{m}b_{n}^{(i)}\partial^{M+N-m-n-i}\\
  \hspace{3cm}=\sum_{l=0}^{\infty}(\sum_{m=0}^{l}\sum_{i=0}^{l-m}
{M-m \choose i} a_{m}b_{l-m-i}^{(i)})\partial^{M+N-l}.
\end{array}
$$
\begin{definition}
$E$ is called the {\em ring of formal pseudo-differential operators with 
coefficients in $\power{R}{x}$}. 
\end{definition}
{\bf Remark 1}\freespace
In our notation we follow M.~Mulase \cite{M1}. Other authors use the name 
of {\em micro-differential operators} for the objects which we call 
formal pseudo-differential operators.
\freeline\\
{\bf Remark 2}
\begin{itemize}
\item
{}From the above construction it is clear that $E$ is an associative, 
non-commutative ring, which has the additional
 structure of a left $\power{R}{x}$-module. 
\item
$E$ has a filtration by left 
$\power{R}{x}$-submodules
$$E^{(m)}:=\left\{\sum_{n\in \Bbb{Z}} f_{n}\partial^{n}/ f_{n}\in \power{R}{x}, 
f_{n} = 0 \textrm{ for } n>m\right\}.$$
The {\em order} of an element $P\in E$ is defined to be the minimum 
$m\in \Bbb{Z}$ such that $P\in E^{(m)}$. In particular, the order of an 
element of $D$ coincides with its degree
 when we consider it as a polynomial in the variable $\partial$.

For an operator $P=\sum_{n=0}^{\infty} f_{n}\partial^{N-n}$ of order $N$,
$f_{0}$ is called its {\em leading coefficient}.
\item
An operator $P\in E$ can be written in the {\em right normal form} 
$P= \sum_{m=0}^{\infty}a_{m}\partial^{M-m}$ or in the {\em left normal form} 
$P= \sum_{n=0}^{\infty}\partial^{N-n}b_{n}$. 
 It is an easy consequence of the Leibniz rule that the order of an 
operator is the same in the left and the right normal form and the 
leading coefficient does not change. So we see that $E$ is also a right 
$\power{R}{x}$-module and that $E^{(m)}$ 
gives rise to a filtration of $E$ by right $\power{R}{x}$-submodules.

For more properties of 
formal pseudo-differential operators see the appendix \ref{B}.
\end{itemize}

$E$ contains the right ideal $xE$ generated by $x$. Denote by $\sigma: 
E\rightarrow E/xE$ the projection. So we have for formal pseudo-differential 
operators in the right normal form:
$$
\sigma(\sum_{n=0}^{\infty}  f_{n} \partial^{N-n}) =
\sum_{n=0}^{\infty} f_{n}(0) \partial^{N-n}\in 
\negpower{R}{\partial^{-1}} = E/xE.$$
Let us set
$$y := \partial^{-1}.$$
Obviously we get, for all $n\in\Bbb{Z}$:
$$\sigma(E^{(n)}) = \power{R}{y}\cdot y^{-n}.$$
\begin{definition}
\label{act PDO}
The projection map $\sigma$ defines an action of $E$ on $\negpower{R}{y}$ 
as follows:\\
Take $P\in E$ and $v\in\negpower{R}{y}$. Then there is an operator $Q\in E$ 
such that $v=\sigma(Q)$. Define
$$P(v) := \sigma(QP).$$
\end{definition}
{\bf Remark}
\begin{itemize}
\item This definition does not depend on the choice of $Q$. Note that for an 
invertible operator $P\in E$, $P^{-1}:\negpower{R}{y} \rightarrow 
\negpower{R}{y}$ is inverse to the map $P:\negpower{R}{y} \rightarrow 
\negpower{R}{y}$. If $P$ is an operator of order 0 with invertible 
leading coefficient, then $P$ is invertible and the
induced map is an automorphism preserving orders, i.e., for all $n\in\Bbb{Z}$,
$$ P : \power{R}{y}\cdot y^{n} \stackrel{\sim}{\rightarrow} 
\power{R}{y}\cdot y^{n}.$$
\item If $P$ is a formal pseudo-differential operator with constant coefficients, 
$P\in \negpower{R}{\partial^{-1}}$, 
then we can regard $P$ as an element of $\negpower{R}{y}$ and it is easy to see 
that in this case the action of $P$
 on $\negpower{R}{y}$ coincides with the usual multiplication in $\negpower{R}{y}$:
$$\sigma(QP) = \sigma(Q)\cdot P= P\cdot \sigma(Q).$$
\item At this point, our way differs slightly from the one taken by M.~Mulase 
\cite{M1}. There, the quotient is taken by
 $Ex$ and operators $P\in E$ act on $\negpower{R}{y}$ from the left. Our 
approach will find its justification in 
  Section \ref{eigenvalue}.
\end{itemize}
\begin{proposition}
Let $P\in E$ be an operator of order 0 with invertible leading
coefficient and
$P:\negpower{R}{y}
\rightarrow \negpower{R}{y}$ the induced automorphism defined above. 
Then $P$ induces an automorphism
$$P:\frak{G}_{F,\alpha}^{1}(Spec(R)) \rightarrow 
\frak{G}_{F,\alpha}^{1}(Spec(R))$$
for all integers  $\alpha\in \Bbb{Z}$ and all elements $F\in K(Spec(R))$.
\end{proposition}
{\bf Proof}\freespace
cf. \cite{M1}, Prop. 4.2.
\begin{definition}
A formal pseudo-differential operator $T\in E$ is called {\em admissible} if 
it is an operator of order 0 with invertible leading coefficient such that
$$T\partial T^{-1} \in \negpower{R}{\partial^{-1}}.$$
The group of admissible operators is denoted by $\Gamma_{a}$.
\end{definition}
\begin{lemma}
\label{admissible}
\begin{enumerate}
\item An operator $T$ is admissible if and only if it has the form 
$$T = exp(c_{1}x)\cdot(\sum_{i=0}^{\infty} f_{i}\partial^{-i}),$$
where $c_{1}\in R$, $f_{i}\in\power{R}{x}$ is a polynomial of degree at 
most $i$, and $f_{0}\in R$ is invertible.
\item Let $v\in\negpower{R}{\partial^{-1}}$ be a monic element of order 
$-r$, $r\neq 0$. Then there is an admissible operator $T\in\Gamma_{a}$ such 
that 
$T\partial^{-r}T^{-1} = v$.
\end{enumerate}
\end{lemma}
{\bf Proof}\freespace
This is a direct consequence of Lemma \ref{XLX} and its corollaries.
\begin{definition}
Two embedded Schur pairs $(A_{1},W_{1})$ and $(A_{2},W_{2})$ of rank $r$, 
index $F$ and level $\alpha$ are said to be {\em equivalent} if there 
is an admissible operator T such that
$$T^{-1}A_{2} T = A_{1}, \quad TW_{2} = W_{1},$$
where  $A_{1}$ and $A_{2}$ are understood to be  subalgebras of 
$\negpower{R}{\partial^{-1}}$, i.e., $T$ acts by conjugation,
 while $W_{1}$ and $W_{2}$ are understood to be subspaces of $\negpower{R}{y}$ 
and  the action of $T$ on $W_{2}$ is defined by Definition \ref{act PDO}. 
So we get $\frak{E}_{\alpha}\frak{S}'^{r}_{F}(Spec(R))/\Gamma_{a}$.
\end{definition}

Now let $(A,W)$ be an arbitrary embedded Schur pair of rank $r$. By Lemma 
\ref{pure}, $A\subset\negpower{R}{z}$ for some element $z\in\negpower{R}{y}$ 
of order $-r$. By Lemma \ref{admissible} there is an operator 
$T\in\Gamma_{a}$ such that
$$ T^{-1}zT = y^{r}.$$
Consequently, $A\subset\negpower{R}{z}$ implies that
$$T^{-1}AT \subset T^{-1}\negpower{R}{z}T = \negpower{R}{T^{-1}zT}=
\negpower{R}{y^{r}}.$$
This proves
\begin{lemma}
Every equivalence class of embedded Schur pairs contains a representative 
which corresponds to a Schur pair. \QED
\end{lemma}
\begin{definition}
Let $\alpha\in\Bbb{Z}$ be an integer. Two Schur pairs of rank $r$ and index $F$ 
are said to be {\em $\alpha$-equivalent} if the associated embedded Schur 
pairs of level $\alpha$ are equivalent.

We call two geometric data {\em $\alpha$-equivalent} if the corresponding 
Schur pairs are $\alpha$-equivalent.
\end{definition}
{\bf Remark}\freespace
Observe that the $\alpha$-equivalence depends on the congruence class of 
$\alpha$ modulo $r$.
\freeline

Now let us study  the $(-1)$-equivalence in more detail.
\begin{lemma}
Let $T\in\Gamma_{a}$ be such that for two algebras $A_{1}$ and $A_{2}$ of 
pure rank $r$ which are both contained in  
$\negpower{R}{y^{r}}$:
$$T A_{1} T^{-1} = A_{2}.$$
Then $T \power{R}{y^{r}} T^{-1} = \power{R}{y^{r}}$.

Conversely, for every algebra $A_{1}$ of pure rank $r$ with $A_{1}
\subseteq \negpower{R}{y^{r}}$ and every $T\in\Gamma_{a}$ satisfying 
$T \power{R}{y^{r}} T^{-1} = \power{R}{y^{r}}$, $TA_{1}T^{-1} 
\subseteq \negpower{R}{y^{r}}$ is also an algebra of pure rank $r$.
\end{lemma}
{\bf Proof}\freespace
The last part is obvious. As for the first one, remember that 
there are elements $a,b\in A_{1}$ such that
$$a^{-1}b = y^{r} + \sum_{i\geq 2} \alpha_{i}y^{ir}$$
(cf. Lemma \ref{pure}). Then $TA_{1}T^{-1} = A_{2} \subseteq 
\negpower{R}{y^{r}}$ implies:
$$
\begin{array}{rcl}
Ta^{-1}bT^{-1} & = & Ta^{-1}T^{-1}TbT^{-1}\\
&=& (TaT^{-1})^{-1} TbT^{-1}\\
& \in & \power{R}{y^{r}}
\end{array}
$$
hence $T\power{R}{a^{-1}b} T^{-1} \subseteq \power{R}{y^{r}}$. 
But $\power{R}{y^{r}}$ equals  $\power{R}{a^{-1}b}$, and so we obtain $T 
\power{R}{y^{r}} T^{-1} \subseteq \power{R}{y^{r}}$. The second 
inclusion we get using the fact that 
$A_{1}= T^{-1}A_{2}T$.
\QED
\freeline

Now let us consider the action of $T\in\Gamma_{a}$ on the second 
component of an embedded Schur pair $(A,W)$ of level $-1$
$$W\subset \negpower{R}{y} = \bigoplus_{i=1}^{r} \negpower{R}
{y^{r}}\cdot y^{i}.$$
For $c_{i}\in \negpower{R}{y^{r}}$ we have
$$T(\sum_{i=1}^{r} c_{i} y^{i}) = 
\sum_{i=1}^{r} (T^{-1}c_{i} T)T(y^{i}).
$$
In particular, by Corollary \ref{XLX2}, $T$ is determined by $A$ and $T^{-1}AT$ 
up to an operator with constant 
coefficients. So we obtain:
\begin{proposition}
\begin{enumerate}
\item
If two geometric data $(C,\pi,Spec(R),P,\rho_{1},{\cal F},\Phi_{1})$ and 
$(C,\pi,Spec(R),$ $P,\rho_{2},{\cal F},\Phi_{2})$ of rank $r$ and 
index $F$ are $(-1)$-equivalent then 
there is an automorphism of $R$-algebras
$$h: \power{R}{z} \stackrel{\sim}{\rightarrow}  \power{R}{z}$$
satisfying $\rho_{2} = h\circ \rho_{1}$. 
\item Two geometric data
 $$(C,\pi,Spec(R),P,\rho,{\cal F},\Phi_{1}) \textrm{ and }
 (C,\pi,Spec(R),P,\rho,{\cal F},\Phi_{2})$$
 of rank $r$ and 
index $F$ are $(-1)$-equivalent if and only if there are elements
 $d_{1},\ldots,d_{r-1}\in \power{R}{z}$ and 
$d_{0}\in \power{R}{z}^{*}$ such that for 
$$
M = \left(
\begin{array}{*{5}{c}}
d_{0} & d_{1}&\ldots &d_{r-2}&d_{r-1}\\
d_{r-1}y^{r}& d_{0} & \ldots & d_{r-3}& d_{r-2}\\
\multicolumn{5}{c}{\dotfill}\\
d_{1}y^{r}& d_{2}y^{r}&\ldots& d_{r-1}y^{r}& d_{0}
\end{array}
\right)
$$
$$\rho\Phi_{2} = M\circ\rho \Phi_{1}.$$

\end{enumerate}
\QED
\end{proposition}

\begin{corollary}
Note that in the case $r=1$,  two geometric data $(C,\pi,Spec(R),P,\rho_{1},
{\cal F},\Phi_{1})$ and $(C,\pi,Spec(R),P,\rho_{2},{\cal F},\Phi_{2})$ 
are always $(-1)$-equivalent.
\QED
\end{corollary}

\subsection{Classification of commutative algebras of differential operators}
After these preliminaries we now come to the main object: the 
classification of commutative algebras of ordinary differential 
operators with coefficients in $\power{R}{x}$, where $R$ is a 
commutative noetherian $k$-algebra, for some field $k$ of 
characteristic zero. 

\begin{definition}
A commutative subalgebra $B$ of $D$ is said to be {\em elliptic of 
pure rank $r$} if
\begin{itemize}
\item $r=gcd(ord(P)/P\in B)$;
\item There are monic elements $P,Q\in B$ such that $gcd(ord(P),ord(Q))=r$.
\end{itemize}
The set of all such subalgebras of $D$ is denoted by ${\cal B}_{r}(R)$.
\end{definition}
Let us start  the observations with the following
\begin{lemma}
\label{XB}
If $B$ is an elliptic subalgebra of rank $r$ of $D$ then there is a 
 formal pseudo-differential operator $X$ of order 0 with coefficients in 
$\power{R}{x}$ and invertible leading coefficient such that
$$A:=X^{-1} B X \subseteq \negpower{R}{\partial^{-1}} $$
and $A$ is an algebra of pure rank $r$.
\end{lemma}
{\bf Proof}\freespace
Choose a monic operator $P\in B$ of order $N$ greater than 0. Then by Lemma 
\ref{XLX} there exists a formal pseudo-differential operator $X$  of order 0 
with invertible leading coefficient such
that
$$X^{-1} P X = \partial^{N}.$$
Let $Q\in B$ be an arbitrary element of $B$. Since $P$ and $Q$ commute, we get
$$
\begin{array}{rcl}
0&=& X^{-1}(PQ-QP)X\\
&=& (X^{-1}PX)(X^{-1}QX) - (X^{-1}QX)(X^{-1}PX)\\
&=& \partial^{N}(X^{-1}QX) - (X^{-1}QX)\partial^{N}.
\end{array}
$$
{}From the proof of Corollary \ref{XLX2} one gets that then $X^{-1}QX$ 
must have constant coefficients, i.e., $X^{-1}QX\in 
\negpower{R}{\partial^{-1}}$. Finally, observe that the rank and the 
monicity of an operator are preserved under  conjugation by  
$X$. So, $A$ is in fact an algebra of pure rank $r$.
\QED
\freeline

Let us look how far the name ``elliptic'' is justified for such an 
algebra $B$. 
 Let $\frak{m}$ be a maximal ideal of $R$. Then $B/\frak{m}$  is a 
commutative algebra of differential operators with coefficients in 
$\power{K}{x}$ for some field $K$ containing $k$. A monic operator $P\in B$ of 
positive order, which 
exists by definition, gives us a monic operator $P\in B/\frak{m}$ 
of positive order. From Lemma \ref{XB}, applied to $B/\frak{m}$, we 
know that every $Q\in B/\frak{m}$ must have constant leading coefficient. 
In particular, since $K$ is a field, the leading coefficient of $Q$ is an
 element of $K^{*}$. This implies that $Q$ is really an elliptic ordinary 
differential operator, i.e., $B$ is a family of algebras of elliptic operators 
parametrized by $R$.

However, this is not the only way to interpret $B$. Let us take a 
commutative subalgebra 
$R\subset \multpower{k}{t}{m}\/[\/\frac{d}{dt_{1}}, \ldots , 
\frac{d}{dt_{m}}\/]$. Then $Q\in B$ is a partial differential operator, 
which does {\bf not} need to be elliptic.
\freeline

Let us continue with our construction.

\begin{lemma}[Sato]
\label{Sato1}
A formal pseudo-differential  operator $P\in E$ is a differential operator if 
and only if it preserves $\sigma(D)$ in $\negpower{R}{y}$, i.e.,
$$P\sigma(D) \subseteq \sigma(D).$$
\end{lemma}
In the proof one may follow  \cite{M1}, Lemma 7.2.
\QED
\freeline

Now let us have a look at another definition of infinite Grassmannians
\begin{definition}
The {\em Sato Grassmannian} is defined to be 
$$SG^{+} := \{J\subset E \textrm{closed subspace }/
J\oplus E^{(-1)} = E, DJ\subseteq J\}$$
\end{definition}
This is a relative version of the Grassmannian originally used by Sato. 
There is the following connection to the Grassmannians we have considered 
until now:
\begin{theorem}[Sato]
\label{Sato2}
\begin{enumerate}
\item Let $\Gamma_{m}$ be the group of monic formal pseudo-differential 
operators of order zero and let $SG^{+}$ be the Sato Grassmannian 
defined as above. Then there is a natural bijection $\alpha:\Gamma_{m} 
\stackrel{\sim}{\rightarrow} 
SG^{+}$ given by
$$\Gamma_{m} \ni X \stackrel{\alpha}{\mapsto} \alpha(X) = J = DX^{-1} 
\in SG^{+}.$$
\item Set $\frak{G}^{+}(Spec(R)):=\{W\in \frak{G}^{1}_{0,-1}
(Spec(R))/W\oplus \power{R}{y}y =\negpower{R}{y}\}$. Then the 
projection $\sigma : E\rightarrow \negpower{R}{y}$ induces a bijection 
$$\sigma : SG^{+} \stackrel{\sim}{\rightarrow} \frak{G}^{+}(Spec(R)).$$
\end{enumerate}
\end{theorem}
{\bf Proof}\freespace
Lemma \ref{inv pseudo} implies that $\Gamma_{m}$ is a group. Now we can apply
 the proof of  \cite{M1}, Thm. 7.4., 
and obtain our result.
\QED
\freeline\\
{\bf Remark}\freespace
$\frak{G}^{+}(Spec(R))$ is the generalization of the 
{\em big cell of the Grassmannian of rank 1, index 
0 and level $-1$} (cf. \cite{M1}).
\begin{definition}
We call two algebras $B_{1},B_{2}\in {\cal B}_{r}(R)$ {\em equivalent} 
if there is an invertible element $f\in\power{R}{x}$ such that
$$B_{1}= f B_{2} f^{-1}.$$
We denote by $\bar{{\cal B}}_{r}(R)$ the set of these equivalence classes.
\end{definition}
\begin{theorem}
For all $r\geq 1$, there is a canonical bijection 
$$\mu_{r} : \bar{{\cal B}}_{r}(R) \rightarrow 
\frak{E}_{-1}\frak{S}'^{r,+}_{0}(Spec(R))/\Gamma_{a},$$
where $\frak{E}_{-1}\frak{S}'^{r,+}_{0}
(Spec(R))$ denotes the subset of $\frak{E}_{-1}\frak{S}'^{r}_{0}
(Spec(R))$ consisting of embedded Schur pairs $(A,W)$ with $W\in 
\frak{G}^{+}(Spec(R))$.
\end{theorem}
{\bf Proof}\freespace
Given $B\in {\cal B}_{r}(R)$. Using Lemma \ref{XB} and 
Theorem \ref{Sato2}, we construct $A:=X^{-1}BX\subseteq \negpower{R}
{\partial^{-1}}$ and $W:= X\sigma(D)\in \frak{G}^{+}
(Spec(R))$. $A$ is an algebra of pure rank $r$.

 Since $B$ is contained in $D$, we get $B\sigma(D)\subseteq \sigma(D)$ 
(cf. Lemma \ref{Sato1}). This implies 
$$
A\cdot W  =  X^{-1} B X (\sigma(DX))
= \sigma(DXX^{-1}BX)
= \sigma(DBX)
\subseteq X\sigma(D)
=W.
$$
So we really constructed an embedded Schur pair of the required type. 

If $X_{1}$ is another operator satisfying $X_{1}^{-1}BX_{1}\subseteq 
\negpower{R}{\partial^{-1}}$ then $T:= X^{-1}X_{1}$ is an admissible 
operator and we end up with an equivalent embedded Schur pair. 

Now, what happens if we take $f B f^{-1}$ instead of $B$ for some 
$f\in\power{R}{x}^{*}$? Then: $A'= (fX)^{-1} (f B f^{-1}) (fX) = A$ 
and $W' = (fX)(\sigma(D)) = f_{0}\cdot W = W$, where $f_{0}$ 
denotes the (invertible) constant coefficient of $f$. 
\freeline

As for the inverse way, let us take  $(A,W)\in \frak{E}_{-1}\frak{S}'^{r,+}_{0}
(Spec(R))$. From Theorem \ref{Sato2} we get a unique operator $S\in 
\Gamma_{m}$ such that
 $W = S\sigma(D)$.  Let us define $B:= SAS^{-1} \subseteq E$. 
Using Lemma \ref{Sato1} we obtain that $B\subseteq D$. We only have to check that 
$\sigma(DB)$ is contained in
 $\sigma(D)$, or, equivalently, that 
$S\sigma(DB)$ is a subset of  $S\sigma(D)$:
$$
S\sigma(DB) = \sigma(DBS) = \sigma(DSS^{-1}BS) = A\cdot W \subseteq W = S\sigma(D).$$

If we start with an equivalent embedded Schur pair $T(A,W)$, for some $T\in 
\Gamma_{a}$, we get an operator $S_{1}\in \Gamma_{m}$ such that
$$S_{1} (\sigma(D)) = TW = TS (\sigma(D)).$$
Again Lemma \ref{Sato1} implies  $S_{1}^{-1}TS\in D$. But this 
operator is an invertible formal pseudo-differential operator of order 0. So  
$S_{1}^{-1}TS$ is in fact an invertible element $f$ of $\power{R}{x}$
 and we conclude that
$$B'= S_{1} TAT^{-1} S_{1}^{-1} = 
f S A S^{-1} f^{-1} = fBf^{-1}.$$
So we in fact end up with an equivalent algebra.
\QED
\freeline

We have established a one-to-one correspondence of elliptic commutative 
algebras of differential operators with coefficients in $\power{R}{x}$ 
and certain equivalence classes of embedded Schur pairs. Using the 
results of the sections \ref{affine base} and \ref{pseudo DO} we can now state:
\begin{corollary}
There is a bijection between equivalence classes of commutative elliptic 
subalgebras of $\poly{\power{R}{x}}{\frac{d}{dx}}$ of pure rank $r$ 
 and $(-1)$-equivalence classes of  geometric data  
$$(C,\pi,Spec(R),P,\rho,{\cal F},\Phi)$$ 
of rank $r$ and index $0$ with 
the extra-condition that 
$$H^{0}({\cal F}) = H^{1}({\cal F}) = 0.$$
In particular, every  sheaf corresponding to a commutative 
algebra of differential operators is strongly semistable with respect to 
$P$ (cf. Section \ref{geom. properties}).
\QED
\end{corollary}

\subsection{Eigenvalue problems}
\label{eigenvalue}
As a motivation, let us approach the above constructed relation
 from another side. We take an elliptic
 commutative algebra $B$ of differential operators with coefficients in 
$\power{R}{x}$. It is an interesting 
 problem to find out the common eigenfunctions of all operators belonging 
to $B$. Let $f\in\power{R}{x}$ be such 
 a common eigenfunction, i.e., for each $P\in B$, there is some 
$\lambda(P)\in R$ such that
$$P(f) = \lambda(P)\cdot f.$$
One easily sees that, for a given $f$, the map
$$
\begin{array}{ccccc}
\lambda & : & B & \rightarrow & R\\
&& P & \mapsto & \lambda(P)
\end{array}
$$
is a homomorphism of $R$-algebras.

On the other hand, given a homomorphism $\lambda  :  B  \rightarrow  R$, what 
are the eigenfunctions of $B$ with
 respect to $\lambda$?

Let $(A,W)$ be an embedded Schur pair corresponding to $B$, 
$W=\sigma(DS)$, $A=S^{-1}BS$, for a formal
 pseudo-differential operator $S$ of order zero with invertible leading
coefficient. For
a
given function $f\in\power{R}{x}$, we define an $R$-linear map \begin{equation}
\label{eigenf}
f  :  W  \rightarrow  R
\end{equation}
by 
$$f(\sigma(QS)) := \sigma(Q(f)), \textrm{ for } Q\in D.$$
This map is well-defined. For, if $\sigma(QS) = \sigma(Q'S)$, then $\sigma(Q) = 
\sigma(Q')$, i.e., $Q-Q'\in xE$. Consequently, 
$\sigma((Q-Q')(f)) = 0$. 

Now we claim
\begin{proposition}
$f\in\power{R}{x}$ is a common eigenfunction of the elements of $B$ with 
the eigenvalue $\lambda$ if and only if $\lambda$
 makes the map (\ref{eigenf}) $A$-linear, i.e., for $a\in A$, $a=S^{-1}PS$, $P\in B$: 
$$f(a\cdot w) = \lambda(P)\cdot f(w).$$
\end{proposition}
{\bf Proof}\freespace
First assume that $f$ is an eigenfunction as above. Then $P(f)=\lambda(P)\cdot f$
 for all $P\in B$. Therefore, for
 $w= \sigma(QS)\in W$ and $a=S^{-1}PS\in A$
$$
\begin{array}{rcl}
f(a\cdot w) & = & f(\sigma(QSS^{-1}PS))\\
&=& f(\sigma(QPS))\\
&=& \sigma((QP)(f))\\
&=& \sigma(Q(\lambda(P)\cdot f))\\
&=& \lambda(P)\cdot \sigma(Q(f))\\
&=& \lambda(P)\cdot f(w).
\end{array}
$$
On the other hand, assume that the $A$-linearity holds. Then, for all ordinary 
differential operators $Q$, 
$$\sigma(Q(P(f))) =
\sigma(Q(\lambda(P)\cdot f)).$$
We claim that this implies $P(f) = \lambda(P)\cdot f$. But this is clear, 
since, for every 
$g=\sum_{n\geq 0}g_{n}x^{n}\in\power{R}{x}$, 
$\sigma(\partial^{n}(g)) = n!\cdot g_{n}$, i.e., $g$ is in fact determined 
by $\sigma(Q(g))$, $Q\in D$. 
\QED

\subsection{Examples}
To get a better idea of the formal correspondence 
constructed in this chapter, let us give some easy examples. 
Certainly, the results we are going to obtain are not new. 
They only serve to illustrate the construction.

We start with the easiest case: Let $B\in {\cal B}_{1}(R)$ be 
an algebra containing a monic element of order 1. We prove
\begin{lemma}
All algebras $B\in {\cal B}_{1}(R)$ containing a monic element 
of order 1 are equivalent.
\end{lemma}
{\bf Proof}\freespace
Take $u = \sum_{i\geq 0} v_{i} x^{i} \in \power{R}{x}$. It suffices  
to prove that there is an invertible formal power series $f = 
\sum_{i\geq 0} f_{i} x^{i} \in \power{R}{x}$ such that
$$f^{-1} \partial f = \partial+u.$$
Let us construct $f$: 
$$f^{-1} \partial f = f^{-1} f \partial + f^{-1} f' = \partial + 
f^{-1} f';
$$
therefore we only need to construct $f$ such that $u = f^{-1} f'$, i.e., 
$fu=f'$. We write this expression as a formal power series
$$\sum_{i\geq 0} (\sum_{j=0}^{i} f_{j} u_{i-j}) x^{i} = 
\sum_{i\geq 0} (i+1) f_{i+1} x^{i},$$
so that we get as a necessary and sufficient condition:
$$f_{i+1} = \frac{1}{i+1} \sum_{j=0}^{i} f_{j} u_{i-j},$$
which is solvable for an arbitrarily given $u$.
\QED
\freeline

What does this say in terms of curves and sheaves? 
The curve associated to such a $B$ is obviously  $\Bbb{P}^{1}_{R}$ with 
some section $P$. One possible line bundle with vanishing cohomologies  
is ${\cal O}_{\Bbb{P}^{1}_{R}}(-P)$. The different local parametrizations 
of $\Bbb{P}^{1}_{R}$ and local trivializations of 
${\cal O}_{\Bbb{P}^{1}_{R}}(-P)$ along $P$ factor out under the 
equivalence relation. Now the above lemma reads as:
\begin{corollary}
For any two sections of $\Bbb{P}^{1}_{R}$ there is an isomorphism of 
$\Bbb{P}^{1}_{R}$ mapping  one section into the other. Furthermore, 
given a section $P$, ${\cal O}_{\Bbb{P}^{1}_{R}}(-P)$ is the only  
coherent sheaf  on $\Bbb{P}^{1}_{R}$ of rank 1 with vanishing cohomology groups 
which satisfies the conditions of Definition \ref{def data}.
\QED
\end{corollary}
\freeline

The next interesting case is that of algebras $B\in {\cal B}_{1}(R)$ 
containing monic elements of order 2 and 3, but without any element of 
order 1. These correspond to families of reduced and irreducible curves 
with arithmetic genus 1.  

First take $R=k$. Then we get a single curve. One may ask how the singularity 
of the curve is displayed in the associated algebra of differential operators: 
\begin{proposition}
An algebra $B$ as above corresponds to a singular plane cubic with its point at 
infinity as a section if and only if there is a 
formal pseudo-differential operator $T$ of order 0 
with invertible leading coefficient such that  $TBT^{-1}$ is an algebra of
{\bf differential} operators with constant coefficients.
\end{proposition}
{\bf Remark}\freespace
Keep in mind that this does not say that $B$ is {\em equivalent} to an 
algebra of differential operators with constant coefficients. The 
above transformation changes the sheaf induced by $B$.
\freeline\\
{\bf Proof}\freespace
Let $B$ be given as above. Then for the corresponding Schur pair 
$(A,W)$ we get:
$$A = k\/[\/y^{-2}+\alpha y^{-1}+\beta,  y^{-3} + \gamma y^{-2} + 
\delta y^{-1} + \epsilon\/] $$
with coefficients $\alpha, \beta, \gamma, \delta, \epsilon \in k$. 
Note that of course $\beta$ and $\epsilon$ may be changed arbitrarily, 
and $\gamma$ may be set to 0.

Assume that $\alpha$ is different from 0. We show that such an $A$ 
(understood as an element of ${\cal B}_{1}(k)$ via the identification
 $\partial^{-1}=y$) is in fact equivalent  to an algebra of differential 
operators containing $\partial^{2}$, i.e., we could have chosen $T$ such 
that $TBT^{-1}$ already contains $\partial^{2}$. 

Given $\alpha,  \delta \in k$ we want to construct an invertible power 
series   $$f=\sum_{i\geq 0} f_{i} x^{i} \in \power{k}{x}$$ such that:
\begin{itemize}
\item $f^{-1}\partial^{2}f = \partial^{2}+\alpha \partial+\beta$ for some 
$\beta\in k$;
\item There are numbers $\epsilon, a, b, c \in k$ such that 
$$f^{-1} (\partial^{3} + a \partial^{2} + b \partial + c ) f = 
\partial^{3}  + \delta \partial + \epsilon.$$
\end{itemize}
Let's start the calculation:
$$
f^{-1}\partial^{2}f = f^{-1}(f \partial^{2} + 2f' \partial + f'') = 
\partial^{2} + 2 f^{-1} f' \partial + f^{-1} f''.
$$
We have already seen that the equation $f'=\frac{\alpha}{2}f$ is solvable. 
For  $f$ chosen as such we get:
$$
\begin{array}{l}
f^{-1}\partial^{2}f  =  \partial^{2} + \alpha \partial + 
\frac{\alpha^{2}}{4}  \textrm{  and}\\
\\
f^{-1} (\partial^{3} + a \partial^{2} + b \partial + c ) f=\\ 
 = 
 f^{-1}(f \partial^{3}  + 3 f'\partial^{2} + 3 f''\partial + f'''+ 
af\partial^{2} + 2af' \partial + af'' + bf\partial + 
bf' + cf)\\
= \partial^{3}  + (3 f^{-1}f' + a)\partial^{2} + (3 f^{-1}f'' + 
2a f^{-1} f'+ b)\partial + \\
\hfill + (f^{-1}f'''   + af^{-1}f'' + 
bf^{-1}f' + c)\\
= \partial^{3}  + (\frac{3\alpha}{2} + a)\partial^{2} + 
(\frac{3\alpha^{2}}{4} + 
a\alpha+ b)\partial + (\frac{\alpha^{3}}{8}   + \frac{a\alpha^{2}}{4} + 
 b\frac{\alpha}{2} + c).
\end{array}
$$
{}From this we see that all conditions can be satisfied.
\freeline

We still have to show that the  rings $A = k\/[\/y^{-2}, y^{-3} + 
\delta y^{-1}\/]$ correspond to singular plane cubics and that all 
such cubics occur. 

A singular plane cubic has the equation 
$$z_{0}z_{2}^{2} = z_{1}(z_{1} + \delta z_{0})^{2}$$
with $P=(0:0:1)$ being its point at infinity. 
Note that it has a cusp if $\lambda = 0$ and a node in the 
remaining cases. 
 Using the methods of \ref{Elliptic curves} we get 
$$A = k\/[\/y^{-2}, y^{-3} + \delta y^{-1}\/]$$
as the 
associated ring,
 and the proof is complete.
\QED
\freeline

Now let us have a look at the partial differential equations which are 
produced by the algebras $B$ as considered above. $B$ is generated by 
two elements
$$\begin{array}{rcl}
L&=&\partial^{2} + u \partial + v\\
P&=&\partial^{3} + \alpha \partial^{2} + \beta \partial + \gamma.
\end{array}
$$
We have seen before that we can assume that $u=0$ if we are only interested 
in the equivalence class of $B$. \freeline

 What does it mean for  $B$ to be commutative? $[P,L]=0$ can be interpreted as
$$
\begin{array}{lrcl}
(I)& 2\alpha' & = & 0\\
(II)& \alpha'' + 2\beta' -3v'&=&0\\
(III)& \beta'' + 2\gamma' - 3 v'' - 2\alpha v'&=&0\\
(IV)& \gamma'' - v'''-\alpha v'' -\beta v' & = & 0.
\end{array}
$$
The first equation says that $\alpha$ is a constant. So we can choose 
another normal form 
$$\begin{array}{rcl}
L&=&\partial^{2}  + v\\
P&=&\partial^{3}  + \beta \partial + \gamma,
\end{array}
$$
$v$ and $\gamma$ without constant term. The class of the algebra $B$ uniquely
 determines  $P$ and $L$, and on the other side, these two elements, 
and by (II) and (III) in fact $\beta$, uniquely determine the equivalence 
class of $B$.

Substituting (I), (II) and (III) in (IV) we get
$$\frac{1}{6} \beta''' - \frac{2}{3} \beta\beta' = 0,$$
which is nothing but the stationary Korteveg - de Vries equation. 

We saw that $B$ corresponds to a singular curve with its point at infinity 
if it may be transformed into $\bar{B} = k\/[\/\partial^{2},\partial^{3} + 
\beta \partial \/]$, $\beta$ being a constant. So we get:
\begin{corollary}
A pointed integral curve $(C,p)$ of arithmetic genus 1 is isomorphic to a singular 
plane cubic with its point at infinity if and only if there is a torsion 
free sheaf on $C$   generating a constant solution of the KdV equation. 
The singularity is a node if this constant if different from zero; it is a 
cusp if the constant equals to zero.
\QED
\end{corollary}

After this study of special cases we are now interested in the general 
structure of ${\cal B}_{1}(\Bbb{C})$. 
{}From the correspondence established in the first part of the section 
\ref{Fam. DO} it is clear that two algebras of differential operators 
lead to the same pointed curve if and only if they differ only by 
conjugation with a
formal pseudo-differential operator of order 0 with invertible
leading coefficient. So we obtain:
\begin{proposition}
For any $B\in {\cal B}_{1}(\Bbb{C})$ such that $Spec(B)$ is smooth, 
 the set 
$$\bar{{\cal B}}_{1}(\Bbb{C})(B):=\left\{
\begin{array}{r}
B'\in \bar{{\cal B}}_{1}(\Bbb{C})/ 
\exists\textrm{ pseudo-diff. op. } X 
\textrm{ of order 0 with inv.}\\
\textrm{leading coeff. such that } X^{-1}B'X=B
\end{array}
\right\}
$$
has the structure of a reduced, irreducible, affine 
complex variety of dimension $$g=card\{n\in \Bbb{N}/ \textrm{ there is no 
differential operator of order $n$ in $B$} \}.$$
\end{proposition}
{\bf Proof}\freespace
Let $(C,p)$ be the pointed curve given by $B$. 
 We consider $Pic^{g-1}(C)$. The 
theta divisor in it is given by the condition $h^{0}\neq 0$. 
This divisor is known to be ample, so its complement 
$$U = 
\{ {\cal L} \textrm{ line bundle of degree $g-1$ on } C/ h^{0}
({\cal L}) = h^{1}({\cal L}) = 0 \}$$
 is affine. It is also reduced and irreducible.  Now we take 
the Poincar\'{e} bundle ${\cal P}^{g-1}_{C}$ on $Pic^{g-1}(C) 
\times C$ (normalized with respect to $p$) and restrict it to 
$U\times C$. We choose some local parametrization of $U\times 
C$ near $U\times\{p\}$ and some local trivialization of ${\cal P}^{g-1}_{C}$ along 
$U\times\{p\}$ (for more details see section \ref{Poincare}).
 Different parametrizations and trivializations cancel out under 
the equivalence relation. Now we construct the associated algebra 
\boldmath$B$\unboldmath {} of ordinary differential operators 
with coefficients in $\power{R}{x}$ for 
$U=Spec(R)$. From the universality of the Poincar\'{e} bundle we get 
that  \boldmath$B$\unboldmath {} uniquely parametrizes the 
equivalence classes of algebras 
$B'\in {\cal B}_{1}(\Bbb{C})$ corresponding to the given pointed curve 
$(C,p)$, hence all of $\bar{{\cal B}}_{1}(\Bbb{C})(B)$.
\QED
\freeline\\
{\bf Remark}\freespace
In the case that $g=1$ the so-defined 
\boldmath$B$\unboldmath$\in{\cal B}_{1}(R)$ 
carries a nontrivial family of solutions of the KdV equation.
\appendix
\section{Basic facts on power series rings and formal pseudo-differential 
operators with coefficients in power series rings}
\label{B}

Here we want to sum up some easy properties of power series 
(resp. Laurent series) with coefficients in rings and of formal 
pseudo-differential operators having coefficients 
in rings of power series. 

Assume that $R$ is a commutative ring.
\begin{lemma}
The units of the ring $\power{R}{z}$ are exactly the power series 
with invertible constant coefficient.
\end{lemma}
{\bf Proof}\freespace
Given $s=\sum_{i=0}^{\infty} s_{i}z^{i} \in \power{R}{z}$, we want 
to determine its inverse. Set $t = \sum_{j=0}^{\infty} t_{j}z^{j}$. Then
$$
\begin{array}{rcl}
s\cdot t & = & (\sum_{i=0}^{\infty} s_{i}z^{i})(\sum_{j=0}^{\infty} 
t_{j}z^{j})\\
& = & \sum_{i=0}^{\infty}\sum_{j=0}^{\infty} s_{i} t_{j} z^{i+j}\\
& = & \sum_{l=0}^{\infty}(\sum_{j=0}^{l} s_{l-j} t_{j}) z^{l}.
\end{array}
$$
Hence, $s\cdot t = 1$ if and only if 
\begin{itemize}
\item $1=s_{0}t_{0}$;
\item $\forall l\geq 1,   
s_{0}t_{l} =-\sum_{j=0}^{l-1} s_{l-j} t_{j}$.
\end{itemize}
One can construct an element $t\in\power{R}{z}$ with these properties 
if and only if the constant term $s_{0}$ of $s$ is invertible.
\QED
\begin{lemma}
Let $N$ be an integer and assume that the order 
$M$ of a monic element $s\in \negpower{R}{z}$ is divisible by $N$. Then 
there is a uniquely determined monic element $t\in \negpower{R}{z}$ of 
order $\frac{M}{N}$ such that
$$t^{N} = s.$$
\end{lemma}
{\bf Proof}\freespace
Set $s= z^{M}(1+\sum_{i=1}^{\infty} s_{i}z^{i})$ and  
$t= z^{M/N}(1+\sum_{j=1}^{\infty} t_{j}z^{j})$. Then $t^{N} = s$ if and 
only if
$$
1+\sum_{i=1}^{\infty} s_{i}z^{i} - (1+\sum_{j=1}^{\infty} 
t_{j}z^{j})^{N} = 0.$$
For $n\geq 1$,  the coefficient of $z^{n}$ of this term is 
$$0 = s_{n} - t_{n} - p_{n},$$
where $p_{n}$ is a polynomial in $t_{1},\ldots,t_{n-1}$. Therefore, all 
$t_{j}$ are well-determined by the coefficients $s_{j}$.
\QED
\freeline

{}From now on we assume that $R$ is a $k$-algebra for some field $k$ of
 characteristic zero. 
\begin{lemma}
\label{inv pseudo}
A formal pseudo-differential operator 
$P=\sum_{n=0}^{\infty}s_{n}\partial^{-n}$
of order 0
is invertible by another formal pseudo-differential
operator of order 0 if and only if
its leading coefficient $s_{0}\in\power{R}{x}$ is invertible.
\end{lemma}
{\bf Proof}\freespace
Set $Q:=\sum_{n=0}^{\infty}t_{n}\partial^{-n}$
 and $S:=\sum_{n=0}^{\infty}r_{n}\partial^{-n}$. 
Now let us write down $QP$.
$$
QP = \sum_{l=0}^{\infty}(\sum_{m=0}^{l}\sum_{i=0}^{l-m} {-m \choose i}
t_{m} s_{l-m-i}^{(i)})\partial^{-l}.
$$
If $QP=1$ then $t_{0}s_{0}=1$, i.e., $s_{0}$ must be a unit in 
$\power{R}{x}$. 
Now let us assume that $s_{0}$ is invertible. From the equation above 
one sees that $QP=1$ if and only if
\begin{itemize}
\item $t_{0}=s_{0}^{-1}$;
\item For $l\geq 1$: $t_{l}= s_{0}^{-1}\cdot f_{l}$, where $f_{l}$ is
 a polynomial in $t_{0},\ldots, t_{l-1}$ and $s_{0},\ldots, s_{l}$ and 
 its derivatives. 
\end{itemize}
These equations are solvable, i.e., such coefficients $t_{i}\in 
\power{R}{x}$ exist.

On the other hand, 
$$PS= \sum_{l=0}^{\infty}(\sum_{m=0}^{l}\sum_{i=0}^{l-m} 
{-m \choose i}
s_{m} r_{l-m-i}^{(i)})\partial^{-l}$$
equals  1 if and only if 
\begin{itemize}
\item $r_{0}=s_{0}^{-1}$;
\item For $l\geq 1$: $r_{l}= s_{0}^{-1}\cdot g_{l}$, $g_{l}$ 
being a polynomial in $t_{0},\ldots, t_{l-1}$ and its derivatives 
and $s_{0},\ldots, s_{l}$. 
\end{itemize}
So there are formal pseudo-differential operators $Q$ and $S$ such that 
$QP=1=PS$. From the associativity of $E$ we get $Q=R=P^{-1}$.
\QED
\begin{lemma}
\label{XLX}
Let $L=\sum_{i=0}^{\infty} u_{i}\partial^{N-i}$ be a monic 
formal pseudo-differential operator  of order $N$, $N\neq 0$, with coefficients 
in $\power{R}{x}$. Then there is a formal pseudo-differential operator
$X=\sum_{i=0}^{\infty} s_{i}\partial^{-i}$ of order 0 with
invertible leading coefficient such that
$$X^{-1} L X=\partial^{N}.$$
\end{lemma}
{\bf Proof}\freespace
$X^{-1} L X = \partial^{N}$ is equivalent to $L X = X\partial^{N}$, i.e.,
$$
\begin{array}{rcl}
\sum_{i=0}^{\infty} s_{i}\partial^{N-i} & = & 
(\sum_{i=0}^{\infty} u_{i}\partial^{N-i})
(\sum_{i=0}^{\infty} s_{i}\partial^{-i})\\
&=& \sum_{l=0}^{\infty}(\sum_{m=0}^{l}\sum_{i=0}^{l-m}
{N-m \choose i} u_{m}s_{l-m-i}^{(i)})\partial^{N-l}.
\end{array}
$$
Comparing coefficients, this gives us:
$$
s_{l} = 
\sum_{m=0}^{l}\sum_{i=0}^{l-m}
{N-m \choose i} u_{m}s_{l-m-i}^{(i)},
$$
and so
$$
0 =\sum_{m=0}^{l}\sum_{i=0}^{l-m} {}_{{}_{(m,i)\neq (0,0)}}
{N-m \choose i} u_{m}s_{l-m-i}^{(i)}.
$$
For $(m,i)=(0,1)$ or $(1,0)$, the term $s_{l-1}$ is involved. 
In all other terms, $s$ occurs with lower index. So we get
$$N u_{0}s_{l-1}' + u_{1}s_{l-1} = P_{l-2},$$
where $P_{l-2}$ is a linear combination of the $s_{0},\ldots,s_{l-2}$ and 
its derivatives with coefficients consisting of integer multiples of
 $u_{0},\ldots,u_{l}$.  Using the fact that $u_{0}=1$ we get the equation
\begin{equation}
\label{normation}
s_{l-1}' = -\frac{1}{N} u_{1}s_{l-1} + \frac{1}{N}P_{l-2}.
\end{equation}
Let us write 
$$s_{l-1}=\sum_{i=0}^{\infty}s_{l-1,i} z^{i},
\quad u_{1}=\sum_{i=0}^{\infty}u_{1,i}z^{i},
\quad P_{l-2}= \sum_{i=0}^{\infty}P_{l-2,i}z^{i}$$
with $s_{l-1,i},u_{1,i},P_{l-2,i}\in R$.
 Then the equation (\ref{normation}) can be written as
$$
\sum_{i=0}^{\infty} s_{l-1,i+1}(i+1)z^{i}  = 
-\frac{1}{N}(\sum_{i=0}^{\infty}\sum_{j=0}^{i} u_{1,i-j} s_{l-1,j}z^{i}) 
+ \frac{1}{N}\sum_{i=0}^{\infty} P_{l-2,i}z^{i}.
$$
Equivalently, for all $i=0,1,2,\ldots$:
$$(i+1)s_{l-1,i+1} = -\frac{1}{N}\sum_{j=0}^{i} u_{1,i-j}s_{l-1,j} 
+ \frac{1}{N} P_{l-2,i}.$$
So we see that we can construct an operator $X$ with the required 
properties.
\QED
\begin{corollary}
Assume that $L$ as above is an operator with constant coefficients. 
Then the operator $X$ constructed in the lemma \ref{XLX} has the form
$$X= exp(-\frac{1}{N}u_{1}x)\sum_{i=0}^{\infty} f_{l}\partial^{-l},$$
where $f_{l}$ is a polynomial of degree at most $l$.
\end{corollary}
{\bf Proof}\freespace
We consider again the equation (\ref{normation}). In the given case, it 
has a solution of the form 
$$s_{l-1}= exp(-\frac{1}{N}u_{1}x)f_{l-1}$$
with $exp(-\frac{1}{N}u_{1}x)f'_{l-1}= \frac{1}{N} P_{l-2}$. Set $l=1$:
$$s_{0}=exp(-\frac{1}{N}u_{1}x)c, \quad c\in R^{*}.$$
Inductively, we conclude that 
$$P_{l-2} = exp(-\frac{1}{N}u_{1}x)\cdot g_{l-2}$$
for some polynomial $g_{l-2}$ of degree at most $l-2$. 
So $f_{l-1}$ is in fact a polynomial of degree $\leq l-1$.
\QED
\freeline
\begin{corollary}
\label{XLX2}
The operator $X$ in Lemma \ref{XLX} is determined up to an operator with 
constant coefficients; i.e., if two operators $X_{1}$ and $X_{2}$ satisfy:
$$X_{1}^{-1}LX_{1} = X_{2}^{-1}LX_{2} = \partial^{N}$$
then there is an operator $X$ with constant coefficients such that 
$X_{1}X=X_{2}$.
\end{corollary}
{\bf Proof}\freespace
Denote $X:=X_{1}^{-1}X_{2}$. We show that $X$ has constant 
coefficients. By assumption, $\partial^{N} = X \partial^{N} X^{-1}$, i.e., 
$X\partial^{N}=\partial^{N}X$. So
$$
\begin{array}{rcl}
\sum_{l\geq 0} s_{l} \partial^{N-l} & = & 
\partial^{N}(\sum_{i\geq 0} s_{i} \partial^{-i})\\
&=& \sum_{l\geq 0} \sum_{i=0}^{l} {N \choose i} 
s_{l-i}^{(i)}\partial^{N-l}.
\end{array}
$$
This implies that $\sum_{i=1}^{l} {N \choose i} s_{l-i}^{(i)} = 0$. 
Now we proceed inductively to show: $s_{m}'= 0$ for all $m$:
\begin{itemize}
\item $m=0$: Set $l=1$: ${N \choose 1} s_{0}'=0$, hence $s_{0}'=0$.
\item Assume that $s_{0}'=\ldots=s_{m-1}'=0$. Then of course also 
the higher derivatives of these terms vanish and for $l=m+1$ we again get:
${N \choose 1} s_{m}'=0$.
\end{itemize}
\QED
\freeline

\end{document}